\documentclass[11pt]{amsart}
\usepackage{amsmath}
\usepackage{amssymb}
\usepackage{amsfonts}
\usepackage{amsthm}
\usepackage{verbatim}
\usepackage{stackrel}
\usepackage{cite}
\usepackage{leftidx}
\usepackage{hyperref}
\usepackage{bookman}
\DeclareMathOperator{\Bl}{\mathcal {B}\l}
 
\def\inv{^{-1}}

\DeclareMathOperator{\mgbar}{\overline\M_g}
\DeclareMathOperator{\del}{\partial}
 \DeclareMathOperator{\Hom}{Hom}

\def\refp #1.{(\ref{#1})}

\newcommand{\A}{\mathcal{A}}

\newcommand{\M}{\mathcal{M}}

\newcommand{\Cal}[1]{\mathcal #1}

\def\sch{\text{sch}}
\def\sbr #1.{^{[#1]}}
\def\sfl #1.{^{\lfloor #1\rfloor}}

\def\inv{^{-1}}
\def\?{{\bf{??}}}

\def\dgla{differential graded Lie algebra\ }

\def\Hilb{\text{Hilb}}

\def\Proj{\textrm{Proj}}

\def\M{\Cal M}
\def\A{\Bbb A}

\def\C{\mathbb C}
\def\P{\mathbb P}

\def\ord{\text{\rm ord}}

\def\sym{\text{\rm Sym} }

\def\Spec{\text{\rm Spec} }

\def\ls{\vskip.25in}
\def\ss{\vskip.12in}

\def\O{\mathcal O}

\def\y{\bar{y}}

\def\Sym{\textrm{Sym}}

\def\g{\mathfrak g}

\def\m{\mathfrak m}

\def\hom{\mathfrak {hom}}

\def\1/2{\frac{1}{2}}

\def\I{\mathcal{ I}}

\def\Ann{\textrm{Ann}}

\def\2{{[2]}}
\def\l{\ell}
\def\nl{\newline}

\def\t{\mathcal{T}}

\def\hom{\mathcal{H}\mathit{om}}

\def\<{\langle}
\def\>{\rangle}
\def\sela{SELA\ }

\def\2{{[2]}}
\def\l{\ell}

\def\Proj{\text{Proj}}
\def\scl #1.{^{\lceil#1\rceil}}
\def\spr #1.{^{(#1)}}
\def\sbc #1.{^{\{#1\}}}
\def\supp{\text{supp}}
\def\subpr#1.{_{(#1)}}

\def\beq{\begin{equation*}}
\def\eeq{\end{equation*}}
\newcommand{\td}{\tilde }

\def\g3{{\Gamma\spr 3.}}

\newcommand{\beql}[2]{\begin{equation}\label{#1}#2\end{equation}}

\newcommand{\eqspl}[2]{
\begin{equation}\label{#1}
\begin{split}
#2\end{split}\end{equation}}
\newcommand{\eqsp}[1]{\begin{equation*}
\begin{split}#1\end{split}\end{equation*}}

\newcommand{\exseq}[3]{
0\to #1\to #2\to #3\to 0
}
\newcommand{\beginalphaenum}{
\begin{enumerate}\renewcommand{\labelenumi}{ }
\item \begin{enumerate}
}

\def\eex{\end{rm}\end{example}}
\newcommand\newsection[1]{\section{#1}\setcounter{equation}{0}
}

\newtheorem*{bthm}{Blowup Theorem}
\newtheorem*{nsthm}{Node Scroll Theorem}
\newtheorem{thm}{Theorem}  [section]

\newtheorem*{thm*}{Theorem}
\newtheorem*{prop*}{Proposition}
\newtheorem{cor}[thm]{Corollary}
\newtheorem*{cor*}{Corollary}
\newtheorem*{remark*}{Remark}

\newtheorem{lem}[thm]{Lemma}
\newtheorem{lem*}{Lemma}

\newtheorem*{claim*}{Claim}
\newtheorem{prop}[thm]{Proposition}
\newtheorem{propdef}[thm]{Proposition-definition}
\newtheorem{defn}[thm]{Definition}
\theoremstyle{remark}

\newtheorem{rem}[thm]{Remark}
\newtheorem{crit-rem}[thm]{Critical remark}

\newtheorem{example}[thm]{Example}
\newtheorem*{example*}{Example}

\newtheorem*{defn*}{Definition}

\begin{document}
\title {Structure of the cycle map\\
for Hilbert schemes of families of nodal curves
}

\author
{Ziv Ran}
\date {\today}
\address {\nl UC Math Dept.  \nl
Surge Facility, Big Springs Road\nl
Riverside CA 92521 US\nl ziv.ran @ ucr.edu}
\thanks{\url{http://arxiv.org/0903.3693}}
\begin{abstract}
We study the relative Hilbert scheme
of
 a family of nodal (or smooth)
curves, over a base of arbitrary dimension,
 via its (birational) {\it{ cycle map}}, going to the relative
symmetric product.  We show the cycle map is the blowing up of the
discriminant locus, which consists of cycles with multiple points.
We determine the structure of certain projective bundles called \emph{node scrolls} which play an important
role in the geometry of Hilbert schemes.
\end{abstract}


\maketitle
\begin{footnotesize}
\tableofcontents
\end{footnotesize}
\section*{Introduction}
In the classical (pre-1980) theory of (smooth) algebraic curves, a dominant
role
is played by  divisors-- equivalently,
finite subschemes-- and their parameter spaces, i.e. symmetric
products.
Notably, one of the first proofs of the existence of
special divisors \cite{kl} was based on intersection theory on
symmetric products, developed earlier by Macdonald \cite{macd}. In more
recent developments however, where
the focus has been on moduli spaces of stable curves, subscheme
methods
have been largely absent, replaced by tools related to stable maps
and their moduli spaces (see \cite{Vak} for a sampling stressing
Vakil's work). Our purpose
in this paper (and others in this series) is to develop and apply global
subscheme methods suitable for the study of stable curves and their
families, aiming eventually, inter alia, to extend Macdonald's theory to the
case of families of  curves with at most nodal singularities.
\par
In this paper we always work over the complex numbers.
Fix a family of curves given by a flat projective morphism \beq
\pi:X\to B \eeq over an irreducible  base,
with fibres \beq X_b=\pi\inv(b), b\in B \eeq which are
nonsingular for the generic $b$ and   at worst nodal for every $b$.
For example, $X$ could be the universal family of
automorphism-free curves over the appropriate open subset
of $\overline{\mathcal M}_g$, the moduli space of Deligne-Mumford
stable curves. Consider the
{\it{relative Hilbert scheme}} \beq X\sbr m._B=\Hilb_m(X/B), \eeq
which parametrizes length-$m$ subschemes of $X$ contained in fibres
of $\pi$. This comes endowed with a \emph{cycle map} (also called
 'Hilb-to-Chow'-- in this case, 'Hilb-to-Sym'-- map)
to the relative symmetric product
 \beq c_m:X\sbr
m._B\to X\spr m._B. \eeq
See
 \S \ref{cyclemap-review-sec} for a review.
Because $X\spr m._B$ may be considered 'elementary' (though it's
highly
singular- see \cite{geonodal}), $c_m$ is a natural tool for studying $X\sbr
m._B$.
The structure of $c_m$ is the object of this paper.
Our first main result is the following theorem which was announced with a sketch of proof in \cite{R2},
where some applications are given as well.
\begin{bthm}
$c_m$ is equivalent to the
blowing up of the discriminant locus \beq D^m\subset X\spr
m._B, \eeq which is the Weil divisor parametrizing nonreduced
cycles.
\end{bthm}
  In particular, we obtain an effective
 Cartier divisor
\[2\Gamma\spr m.=c_m\inv(D^m)\] so that
$-2\Gamma\spr m.$ can be identified with the natural $\O(1)$
polarization of the blowup. In fact, we shall see that $\Gamma\spr
m.$ also exists as a Cartier divisor, not necessarily effective, and
the dual of the associated line bundle, i.e. $\O(-\Gamma\spr m.)$,
will be (abusively) called the \emph{discriminant polarization}
(though 'half discriminant' is more accurate); we will also refer to
$\Gamma\spr m.$ itself sometimes the discriminant polarization. We
emphasize that the Blowup Theorem is valid without dimension
restrictions on $B$. As suggested by the Theorem, the discriminant
polarization encodes the additional information in Hilb
vis-a-vis the unwieldy Sym and so, unsurprisingly, plays a central role in
subsequent developments of
geometry and intersection theory on the Hilbert schemes $X\sbr
m._B$.\par The proof of the Blowup Theorem occupies \S\S \ref{preliminary_reductions}-\ref{reverse-engineering-sec} 
and may be outlined as follows.
\par
(i) A preliminary reduction is made to the local case (\S \ref{preliminary_reductions});\par
(iI) we construct an explicit local model $ H$ for the relative Hilbert scheme (\S\ref{local-model-sec});\par 
(ii) we construct an ideal $G$ in the relative Cartesian product, whose
\emph{syzygies} correspond, essentially, to the defining equations of the pullback ${OH}$
of $ H$ over the Cartesian product; this yields
a map $\gamma$ from the blowup of $G$ to ${OH}$ (\S \ref{reverse-engineering-sec});\par
(iv) using the local analysis, it is shown that $\gamma$ is an isomorphism and that $G$ is the ideal
of the ordered discriminant (big diagonal);\par
(v) consequently $\gamma$ descends to an isomorphism from the blowup of the ideal
of the discriminant to $ H$.  \par
 The usefulness of the local model $H$ extends far beyond the
 Blowup Theorem; in particular, it
yields information about the \emph{singularity stratification},
of $X\sbr m._B$, which may be defined as follows. Let $\theta_1,...,\theta_r$ be a
collection of distinct, hence disjoint, relative nodes of the
family, each living in the total space over its own boundary
 component, and let $n_1,...,n_r$ be integers. Set
\[\mathcal S^{m,n.}(\theta_.;X/B)
=\{z:c_m(z)\geq\sum n_i\theta_i\}\subset X\sbr m._B\]
This is mainly interesting when all $n_i\geq 2$. In this case, we construct
a
surjection
\[\bigcup\limits_{1\leq j_i\leq n_i-1,\forall 1\leq i\leq r}
F^{m,n.}_{j.}(\theta.;X/B)\twoheadrightarrow S^{m,n.}(\theta_.;X/B)\]
where each $F^{m,n.}_{j.}(\theta.;X/B)$, called a  \emph{node polyscroll}
(or node scroll, when $r=1$),
is a $(\P^1)^r$-bundle over the smaller Hilbert scheme
$(X^{\theta.})\sbr m-\sum n_i.$, where $X^{\theta.}$ denotes the blowup
(=partial normalization) of $X$ in $\theta_1,...,\theta_r$,
defined over the intersection of the boundary components
corresponding to the $\theta_i$. The fibre parameter of
$i$-th factor of
the node polyscroll encodes a sort of higher-order
(more precisely, $(n_i-1)$-st order) 'slope',
locally at the $i$-th node, and these together constitute the additional
information contained in the Hilbert scheme beyond
what's in the
symmetric product.\par
In the following section 
 \S \ref{flags-sec} we give an analogue of the blowup theorem in the case of flag-Hilbert schemes,
which are often important in inductive arguments and procedures.
\par
Our next main results (see Theorems \ref{nodescroll-thm},
\ref{nodepolyscroll-thm}) determine the structure of node polyscrolls as
$(\P^1)^r$-bundles.
In fact, the disjointness of the nodes (in the total space)
implies that the $\P^1$
factors 'vary independently', which allows us to reduce
to the case of node scrolls, i.e. $r=1$.

Actually, what's essential for the enumerative theory of the Hilbert
scheme,
as studied e.g. in \cite{internodal}, and in which node scrolls play an
essential role,
is the structure of the node
scroll $F$ as a \emph{polarized} $\P^1$ bundle, that is, the  rank-2 vector
bundles $E$ so that there is an isomorphism
$\P(E)\simeq F$, under which the canonical
$\O(1)$ polarization on $\P(E)$
 associated to the
projectivization corresponds to the restriction of the discriminant
polarization $-\Gamma\spr m.$ on $F$. To state the result
(approximately),
denote by $\theta_x, \theta_y$ the node preimages on $X^\theta$,
and by $\psi_x, \psi_y$ the relative cotangent spaces to $X^\theta/T$
along them, and by $[m-n]_*D$, for any divisor $D$ on $X^\theta$,
the 'norm' of $D$, considered as a divisor on $(X^\theta)\sbr m-n.$.
\begin{nsthm}
There is a polarized isomorphism
\[ F^{m,n}_j(\theta)= \P(\O(-D_j^n(\theta))\oplus \O(-D^n_{j+1}(\theta)))\]
where
\[D^n_j(\theta)=-\binom{n-j+1}{2}\psi_x-\binom{j}{2}\psi_y
+(n-j+1)[m-n]_*\theta_x+j[m-n]_*\theta_y+\Gamma\sbr m-n..\]
\end{nsthm}
This result, and its polyscroll analogue,  reduce intersection theory
on polyscrolls to that of the Mumford tautological classes, about which
a great deal is now known thanks to the work of Witten, Kontsevich,
Faber and many others (see e.g. \cite{Vak} and references therein).
The Node Scroll Theorem is one of the main ingredients
of a complete 'Hilbert- tautological' intersection calculus,
 developed in \cite{internodal}, which allows us to extend the intersection
 theory and enumerative geometry
  of a single smooth curve, as developed notably by Macdonald
  \cite{macd} and presented in \cite{acgh}, to the case of families of curves
  with at most nodal singularities, extending work of Cotteril \cite{cotteril1}
  in low degrees. As described in \cite{internodal},
  this intersection calculus has now been implemented on the computer,
  in the form of a Java program called macnodal \cite{macnodal}, 
  due to Gwoho Liu and available from
  the author's web page. See also \cite{grd} for an application to the class of
  the closure of the hyperelliptic class in $\mgbar$.\par
\ss\noindent
{\bf{Acknowledgments}\quad} I thank Mirel Caibar for asking some
stimulating questions early on, and the referee for many helpful, detailed 
comments and suggestions.\nl
{\bf{Convention}}\ \ In this paper we always work over $\C$.
\part{Blowup theorem and discriminant polarization}
\newsection{Review of cycle map}\label{cyclemap-review-sec}
See \cite{ang}, \cite{lehn-montreal} or \cite{sernesi}
for more informatrion.

\subsection{Norms and multisections}
Let $Z=\Spec_T(A)\to T$ be a finite, flat, degree-$m$ morphism  of
algebraic $\C$-schemes, corresponding to a sheaf of $T$- algebras $A$
that is locally $T$-free
of rank $m$.
The action of the algebra $\sym^m_TA$ on the invertible $T$-module
${\bigwedge\limits^m}_T(A)$ yields a $T$-homomorphism of algebras

\[\sym^m_T(A)\to \O_T={\mathrm{End}}_T{(\bigwedge\limits^m}_T(A)).\]
This is a symmetric-tensor version of the norm map, usually given as
a homogeneous polynomial; it can be written locally it terms of determinants.
Applying $\Spec$, we get a $T$-map, called the \emph{canonical
multisection} of $Z/T$,
\[\sigma_{Z/T}:T\to Z\spr m._T=\Spec_T(\sym^m_T(A)).\]
This map is obviously compatible with base-change and satisfies a
'locality' property, namely if $Z=\coprod Z_i$ with each $Z_i$ flat of degree
$m_i$, ,then $\sigma_{Z/T}$ factors through
\[\prod\sigma_{Z_i/T}:T\to \prod Z\spr m_i._T.\]
Consequently, if $t\in T$ and the fibre $Z(t)=\coprod Z_i(t)$ and each
$Z_i$ is supported at a unique point $p_i$, then $\sigma_{Z_i/T}(t)$ is
the unique point of $(Z_i)\spr m_i._T$, usually denoted $m_ip_i$,
and $\sigma_{Z/T}(t)=\sum m_ip_i$.

\subsection{Cycle map}
Let $X\to B$ be a quasi-projective morphism,
$T\to B$ a morphism and $Z$
a $T$-valued point of the relative Hilbert scheme $X\sbr m._B$, i.e. a closed subscheme of
$X\times_B T$ that is finite flat of degree $m$ over $T$.
Examples of possible $T$ include the Hilbert scheme
$X\sbr m._B$ itself and any scheme mapping to it.
We have the canonical multisection, which is a $T$-morphism
\[\sigma_{Z/T}:T\to Z\spr m._T\subset (X\times_BT)\spr m._T=X\spr
m._B\times_BT.\]
Composing with the projection, we get the cycle map, a $B$-morphism
\[c_Z:T\to X\spr m._B.\]
Again, this is compatible with base-change $B'\to B$ and has a locality
property. Moreover, it depends only on $Z$ quasi-intrinsically in the sense that if $Y\subset
X$ is any locally closed subscheme such that $Y\times_BT$ contains $Z$
scheme-theoretically, then $c_Z$ factors through $Y\spr m._B$. Also,
there is an analogous and compatible construction in the analytic
category.

\newsection{Blowup Theorem: Set-up and preliminary
reductions}\label{preliminary_reductions}
\subsection{Set-up} Let
\beq \pi : X\to B \eeq be a flat family of nodal, generically smooth curves
with $X,B$ reduced and irreducible.
Let $X^m_B,  X\spr m._B$, respectively, denote
the $m$th Cartesian and symmetric fibre products of $X$ relative
to $B$. Thus, there is a natural map
\beq \omega_m:X^m_B\to  X\spr m._B \eeq which realizes its target
as the quotient of its source under the permutation action of the
symmetric group $\frak S_n.$ Let \beq \Hilb_m(X/B)=X\sbr m._B \eeq
denote the relative Hilbert scheme parametrizing length-$m$
subschemes of fibres of $\pi$, and
\beq  c=  c_m :X\sbr m._B\to  X\spr m._B \eeq the natural
{{cycle map}} constructed above, associated to the universal
subscheme $Z\subset X\sbr m._B\times_BX$. Let $D^m\subset X\spr
m._B$ denote
the discriminant locus or 'big diagonal', consisting of cycles
supported on $<m$ points (endowed with the reduced scheme
structure). Clearly, $D^m$ is a prime Weil divisor on $ X\spr
m._B$, birational to $X\times_BX\spr m-2._B$ (though it is less
clear what the defining equations of $D^m$ on $ X\spr m._B$ are
near singular points). The main result of Sections 1-4  is the
\begin{thm}[Blowup Theorem]\label{blowup} The cycle map \beq
c_m:X\sbr m._B\to
 X\spr m._B \eeq is equivalent to the blowing up of
$D^m\subset X\spr m._B$.\end{thm}
The proof presented here is an elaboration of the one sketched in \cite{R2}.

\subsection{Reductions}
We begin with some preliminary remarks and reductions.
 To begin with, recall that the cycle map is compatible with base-change,
 as was observed in \S \ref{cyclemap-review-sec}, and note now that the same
 is
 true of the blowup of $D^m$: indeed given a base-change
 $X_{B'}=X\times_BB'$, we have $I_{D^m(X_{B'}/B')}=I_{D^m}\otimes
 \O_{B'}$, hence also
 $I^n_{D^m(X_{B'}/B')}=I^n_{D^m}\otimes \O_{B'}$, so
  \[\bigoplus\limits_n I^n_{D^m(X_{B'}/B')}=(\bigoplus\limits_n
  I_{D^m})\otimes \O_{B'},\] and applying $\Proj$ we get
  \[ \Bl_{D^m}(X\spr m._B)\times_BB'=\Bl_{D^m(X_{B'}/B')}(X_{B'})\spr
  m._{B'}.\]
  \par
 Because the Theorem is local over $B$ and locally any family is a
 base-change from a versal one, we may as well assume $X/B$ is a versal
 deformation of a nodal curve $X_0$, and in particular $X$ and $B$ are
 smooth.\par
 Next, the Theorem is the statement that the natural birational correspondence between $X\sbr m._B$ and $\Bl_{D^m}(X\spr m._B)$
 projects isomorphically both ways (in particular $X\sbr m._B$ is
 irreducible). By GAGA, it suffices to prove for the corresponding
 analytic spaces. Then, since the statement is local over $X\spr m._B$, we may work over a neighborhood of a given cycle $Z=\sum\limits_{i=1}^k m_ip_i$, of
 the form ${\prod}_B(U_i)\spr m_i._B$ where $U_i$ is a suitable analytic
 neighborhood of $p_i$. The corresponding open subset of $X\sbr m._B$ is just ${\prod}_B(U_i\sbr m_i.)_B$, where for an analytic open $U\subset X$,   $U\sbr m._B\subset X\sbr m._B$ is the set of
 schemes contained in $U$. We note that this depends only on $U/B$
 up to analytic isomorphism: e.g. because it can be identified with a Douady space of finite subspaces of $U$; or more directly, by GAGA,
 there is a natural correspondence between analytic families of finite
 subschemes $X/B$ contained in $U$ and finite analytic subspaces of $U/B$. Now choosing $U_i$ appropriately, we may assume there is
 an open subset $V\subset\C^2$ such that $U_i/B$ is a base-change
 of the family $V/T$ given by $xy=t$ (the 'standard model').\par
 Now suppose we could show that $V\sbr m._T, \forall m$, is the
 blowup of $V\spr m._T$ in $D^m$. Then the same is true for $(U_i)\sbr m_i._B, \forall i.$ To conclude that 
 ${\prod} _B (U_i)\sbr m_i._B\simeq
 \Bl_{D^m}{\prod} _B (U_i)\spr m_i._B$, it would suffice to show that
 \[ \Bl_{D^m}{\prod} _B (U_i)\spr m_i._B\simeq  {\prod} _B(\Bl_{D^m} (U_i)\spr m_i._B)\] or equivalently,
 \eqspl{prod-blowup-eq}{ \Bl_{D^m}{\prod} _B (U_i)\spr m_i._B\simeq  {\prod} _{{\prod} _B (U_i)\spr m_i._B}
 (\Bl_{D^m} (U_i)\spr m_i._B)}
 \eqref{prod-blowup-eq} holds because:\begin{enumerate}\item
 The local analysis of the next two sections will show, in particular, 
 that $\Bl_{D^m}V\spr m._T$ is a small blowup, centered over the locus of schemes with multiplicity $\geq 2$ at the node, therefore so is $\Bl_{D^m}(U_i)\spr m_i._B$.
 \item
 \[ I_{D^m(\coprod U_i/B)}|_{{\prod}_B(U_i)\spr m_i._B}
 =\prod I_{D^{m_i}(U_i/B)}.\]
 \item The blowup centers are transverse for different $i$.
 \item The following general remark.

 \end{enumerate}
\begin{rem}\label{product-blowup-rem}Let
$I_1,...,I_k$ be an arbitrary collection of  ideals on a variety $X$, not necessarily 
mutually transverse or even distinct.
\par
 (i) The blowup
$\Bl_{I_1...I_k}X$ of the product ideal is the unique
$X$-dominating component of the fibre product
$\Bl_{I_1}X\times_X...\times_X\Bl_{I_k}X$. For simplicity we check this for
$k=2$. We may work locally over $X$.
If $f_i, g_i$ are generators for $I_1, I_2$ respectively, then the blowup of
$I_1I_2$ is covered by open affines $U_{i,j}$ whose coordinate rings
are generated over $X$ by symbols $[f_{i'}g_{j'}/f_ig_j]$ satisfying
the obvious relations $f_ig_j[f_{i'}g_{j'}/f_ig_j]=f_{i'}g_{j'}$
, $\forall i', j'$.
Similarly with open affines $V^1_i, V^2_j$ for the blowup of $I_1, I_2$,
with generators $[f_{i'}/f_i], [g_{j'}/g_j]$ are regular.
There are obvious maps $U_{i,j}\rightleftarrows V^1_i\times_XV^2_j$,
defined by
$[f_{i'}g_{j'}/f_ig_j]\rightleftarrows [f_{i'}/f_i]\otimes [g_{j'}/g_j] $, leading to
maps over $X$
\[\Bl_{I_1I_2}X\rightleftarrows \Bl_{I_1}X\times_X\Bl_{I_2}X.\]
These clearly give an isomorphism as claimed.\par 
Note that the foregoing argument makes no assumption regarding transversality of $I_1, I_2$.
In general, if $I_1, I_2$
are not transverse, e.g. $I_1=I_2=I$, then $\Bl_{I_1}X\times_X\Bl_{I_2}X$
is reducible: e.g. $[f_1/f_2][f_2/f_1]-1$ is  a zero-divisor
(usually nonzero) on
$\Bl_{I}X\times_X\Bl_{I}X$. The dominating component of
$\Bl_{I}X\times_X\Bl_{I}X$ is
$\Bl_{I^2}X\simeq\Bl_{I}X$.
\par
(ii) In the above situation, if the $\Bl_{I_i}X$ are small blowups,
 i.e. for each $i$ the exceptional locus on $X$ (the center),
 i.e. the non-invertible locus of $I_i$, is of codimension
 $\geq 3$ and its inverse image is of codimension $\geq 2$,
 and if for different $i$ the centers are mutually transverse,
 then the fibre product is in fact irreducible, i.e. has no non-dominating components. This is because any non-dominating component would have to be of smaller dimension, whereas by semi-continuity,
 in the fibre product, which is the inverse image
of the small diagonal in $X^k$ by the
natural map
\[\prod \Bl_{I_i}X\to X^k,\]
every component
is of dimension $\geq\dim(X)$.
\end{rem}
\par
We have now reduced the Theorem to the case where $X/B$ is the standard family $xy=t$, which we assume till further notice; we also
let $U$ denote
any neighborhood of the origin in $X$.
\ls

\newsection{A local model}\label{local-model-sec}
We now give an explicit construction in coordinates
 of the relative Hilbert scheme
of the standard family. This construction will
have many applications beyond the proof of the Blowup Theorem. We
begin with some preliminaries.\subsection{Symmetric product}
Assuming $U/B$ has the local form $xy=t$,
the relative Cartesian product $U^m_B$, as
a subscheme of $U^m\times B$, is given locally  by
\beq x_1y_1=...=x_my_m=t.\eeq Let $\sigma_i^x, \sigma_i^y,
i=0,...,m$ denote the elementary symmetric functions in
$x_1,...,x_m$ and in $y_1,...,y_m$, respectively, where we set
$\sigma_0=1$.
We note that these functions satisfy the relations
\begin{eqnarray}\label{sigma-x-y-rel}
\sigma^y_m\sigma^x_j=t^j\sigma^y_{m-j},\ \
\sigma^x_m\sigma^y_j=t^j\sigma^x_{m-j}, \\
\label{sigma-x-y-rel2}
t^{m-i}\sigma^y_{m-j}=t^{m-i-j}\sigma^x_j\sigma^y_m,\ \
t^{m-i}\sigma^x_{m-j}=t^{m-i-j}\sigma^y_j\sigma^x_m
\end{eqnarray}
(of course the relations in second set follow from those of the first).
Putting the sigma functions together with the projection to $B$, we get
a map
\beq \sigma:U\spr m._B=\Sym^m(U/B)\to \A^{2m}_B= \A^{2m}\times
B\eeq
 \beq \sigma=((-1)^m\sigma_m^x,...,-\sigma_1^x,
(-1)^m\sigma_m^y,...,-\sigma_1^y, \pi^{(m)}) \eeq where $\pi^{(m)}:
X\spr m._B\to B$ is the structure map.

  \begin{lem}\label{sigma-emb} $\sigma$ is an
embedding locally near $mp$ where $p=(0,0)$ is the origin in $U$.\end{lem}\begin{proof} It suffices
to prove this formally, i.e. to show that $\sigma_i^x, \sigma_j^y,
i,j=1,...,m$ generate the completion $\hat{\m}$ of
the maximal ideal of $mp$ in $ X\spr m._B.$ To this end it
suffices to show that any $\frak S_m$-invariant polynomial in the
$x_i, y_j$ is a polynomial in the $\sigma_i^x, \sigma_j^y$ and
$t$. Let us denote by $R$ the averaging or symmetrization operator
with respect to the permutation action of $\frak S_m$, i.e.
\beq R(f)=\frac{1}{m!}\sum\limits_{g\in\frak S_m}g^*(f). \eeq
 Then it suffices to show that the elements
$R(x^Iy^J)$, where $x^I$ (resp. $y^J$) range over all monomials in
$x_1,...,x_m$ (resp. $y_1,...,y_m$) are polynomials in the
$\sigma_i^x, \sigma_j^y$ and $t$.
Because $x_iy_i=t$, we may assume $I, J$ are disjointly supported
in the sense that $I_k>0\Rightarrow J_k=0$.
On the other hand, expanding the product $R(x^I)R(y^J)$
we get a sum of monomials $x^{I'}y^{J'}$ times a rational number;
those with
$I'\cap J'=\emptyset$ add up to $\frac{1}{m!}R(x^Iy^J)$,
while those with $I', J'$ not disjointly supported are divisible by $t$.
Thus,
\beq R(x^Iy^J)-m!R(x^I)R(y^J)=tF \eeq where $F$ is an $\frak
S_m$-invariant polynomial in the $x_i, y_j$ of bidegree
$(|I|-1,|J|-1)$, hence a linear combination of elements of the
form $R(x^{I'}y^{J'}), |I'|=|I|-1, |J'|=|J|-1$. By induction, $F$
is a polynomial in the $\sigma_i^x, \sigma_j^y$ and clearly so is
$R(x^I)R(y^J).$ Hence so is $R(x^Iy^J)$ and we are
done.\end{proof}
\begin{rem} It will follow from the Blowup Theorem \ref{blowup} and its
proof that
the equations (\ref{sigma-x-y-rel}-\ref{sigma-x-y-rel2}) actually
define the image of $\sigma$ scheme-theoretically
(see Cor. \ref{sym-eq} below); we won't need this,
however.\end{rem}
\subsection{A projective family}\label{a projective family}
Now we present a construction of our local model $\tilde H$.
This is motivated by our study in
\cite{Hilb} of the relative Hilbert scheme of a node.
As we saw there, the fibres of the cycle map are
chains consisting of $n$ rational curves where $n$ takes the
values from $n=0$ for the generic fibre (meaning the fibre is a singleton)
to $n=m-1$ for the most special fibre. Therefore, it is reasonable to try to
model the cycle map on a standard pencil  of rational
normal $(m-1)$-tics specializing to a chain of lines. Further
motivation for the construction that follows comes from
\cite{geonodal}, where an explicit construction is given for the
full-flag Hilbert scheme.\par
 Let $C_1,...,C_{m-1}$ be copies of $\P^1$, with homogenous
coordinates $u_i,v_i$ on the $i$-th copy. Let $$\tilde{C}\subset
C_1\times...\times C_{m-1}\times B/B$$ be the subscheme over $B$
defined by
\beql{Ctilde}{ v_1u_2=tu_1v_2,..., v_{m-2}u_{m-1}=tu_{m-2}v_{m-1}.}
This construction is motivated (cf. \cite{geonodal}) by
viewing $u_i/v_i$ as a stand-in for $y_I/x_{I^c}$ where $I\subset [1,m]$
is of cardinality $i$ and $x_I=\prod\limits_{a\in I}x_a$ etc; the
\emph{ratio} is independent of $I$ for fixed $|I|$.
That said,
$\tilde{C}$ is in any event a reduced complete intersection of divisors of
type
\[(1,1,0,...,0), (0,1,1,0,...,0) ,..., (0,...,0,1,1)\]
(relatively over $B$) and it is
easy to check that the fibre of $\tilde{C}$ over $0\in B$ is
\beql{Ctilde0}{ \tilde{C}_0=\bigcup\limits_{i=1}^{m-1}\tilde{C}_i, }
where \beq \tilde{C}_i= [1,0]\times...\times[1,0]\times
C_i\times[0,1]\times...\times[0,1] \eeq and that in a neighborhood of the
special fibre
$\tilde{C}_0$, $\tilde{C}$ is smooth and $\tilde{C}_0$ is its
unique singular fibre over $B.$ We may embed $\tilde{C}$ in
$\P^{m-1}\times B,$ relatively over $B$ using the plurihomogeneous
monomials \beql{Z_i}{ Z_i=u_1\cdots u_{i-1}v_{i}\cdots v_{m-1}, i=1,...,m.
}
These satisfy the relations \beql{Z-rel-quadratic}
{ Z_iZ_j=t^{j-i-1}Z_{i+1}Z_{j-1}, i<j-1} so they embed $\tilde{C}$ as a family
of rational normal
curves $\tilde{C}_t\subset\P^{m-1}, t\neq 0$ specializing to
$\tilde{C}_0$, which is embedded as a nondegenerate, connected
chain
 of $m-1$ lines.\subsection{To Hilb}
 Next consider an affine
space $\A^{2m}$ with coordinates $a_0,...,a_{m-1}$,
$d_0,...,d_{m-1}$. The $a_i, d_j$ are to play the roles of
$\sigma^x_{m-i}, \sigma^y_{m-j}$ respectively
(where as we recall $u_i/v_i$ plays that of
$y_{m-i+1}...y_m/x_1...x_{m-i}$). With this and
the relations \eqref{sigma-x-y-rel}, \eqref{sigma-x-y-rel2} in mind, let
$\tilde{H}\subset\tilde{C}\times\A^{2m}$
be the subscheme defined by
\begin{eqnarray} \label{H-equations} a_0u_1=tv_1, d_0v_{m-1}=tu_{m-1} \\
\nonumber
 a_1u_1=d_{m-1}v_{1},...,a_{m-1}u_{m-1}=d_1v_{m-1}.\end{eqnarray}
 Note that $\tilde H$ comes equipped with a map to $B$ (via the projection to $\Tilde C$), whence
 a projection
 \[p_{\A^{2m}_B}:\tilde H\to A^{2m}_B.\]
Set $L_i=p_{C_i}^*\O(1).$ Then consider the subscheme of
$Y=\tilde{H}\times_{B}U$ defined by the equations
\begin{eqnarray}\label{Y-equations}
F_0:=x^m+a_{m-1}x^{m-1}+...+a_1x+a_0\in
\Gamma(Y,\O_Y)
\\ F_1:=u_1x^{m-1}+u_1a_{m-1}x^{m-2}+...+u_1a_2x+u_1a_1+v_1y
\in\Gamma(Y,L_1) \\ \nonumber ...
\\ \nonumber
F_i:=u_ix^{m-i}+u_ia_{m-1}x^{m-i-1}+...+u_ia_{i+1}x+u_ia_i+
v_id_{m-i+1}y+...+v_id_{m-1}y^{i-1}+ v_iy^i \\
\in\Gamma(Y,L_i)\\ \nonumber ...
\\ F_m:=d_0+d_1y_1+...+d_{m-1}y^{m-1}+y^m\in
\Gamma(Y,\O_Y). \end{eqnarray}

The following statement essentially summarizes results from
\cite{Hilb}.  \begin{thm}\label{Hilb-local}\begin{enumerate}
\item  $\tilde{H}$ is smooth and
irreducible.\item  The ideal sheaf $\I$ generated by
$F_0,...,F_m$ defines a subscheme of $\tilde{H}\times_BX$ that is
flat of length $m$ over $\tilde{H}$ and flat over $X$.
\item The classifying map \beq \Phi=\Phi_\I:
\tilde{H}\to\Hilb_m(U/B) \eeq
is an isomorphism and via $\Phi$, the projection
$p_{\A^{2m}_B}:\tilde H\to \A^{2m}_B$ corresponds to the cycle map.
\item 
$\Phi$ induces an isomorphism
$$\tilde C_0=(\tilde C)_0=p_{\A^{2m}_B}\inv(0)\to
\Hilb^0_m(X_0)=\bigcup\limits_{i=1}^{m-1} C^m_i$$
(cf. \cite{Hilb}) of the fibre of $\tilde{H}$ over $0\in\A^{2m}_B$ with
the punctual Hilbert scheme of the node on the special fibre $X_0$,  in
such a way that the point $[u,v]\in\tilde{C}_i\simeq C_i^m\simeq\P^1$
corresponds to\begin{itemize}\item the subscheme with ideal
$I^m_i(u/v)=(x^{m-i}+(u/v)y^i)\in C^m_i\subset\Hilb^0_m(X_0)$ if $uv\neq
0,$\item the subscheme
$(x^{m+1-i}, y^i)\in C^m_i$ if $[u,v]=[0,1]$,\item  the subscheme $(x^{m-i},
y^{i+1})\in C^m_i$ if
$[u,v]=[1,0].$\end{itemize} In particular, $\tilde C_i$ corresponds
to $C^m_i$.
\item  over the standard open $U_i=(Z_i\neq 0)\subset\P^{m-1}$, a co-basis
    for the universal ideal $\I$ (i.e. a
basis for $\O/\I$) is given by $$1,...,x^{m-i}, y,...,y^{i-1}.$$
\item $\Phi$ induces an isomorphism of the special fibre
$\td H_0$ of $H$ over $B$ with $\Hilb_m(X_0)$, and $\td H_0\subset\td
H$ is a divisor with global normal crossings
$\bigcup\limits_{i=0}^m D^m_i$ where each $D^m_i$ is smooth,
birational to $(x-{\mathrm{axis}})^{m-i}\times (y-\mathrm{axis})^i$, and for
$i=1,...,m-1$ has special
fibre
$C^m_i$ under the cycle map $p_{\A^{2m}_B}$.
\end{enumerate}
\end{thm}\begin{proof} Assertions (i), (ii) are clear
from the defining equations  To prove (iii) and (iv)
consider the point $q_i, i=1,...,m,$ on the special fibre
of $\tilde{H}$ over $\A^{2m}_B$ with coordinates \beq v_j=0,\ \forall
j< i; u_j=0,\ \forall j\geq i. \eeq Then $q_i$ has an affine
neighborhood $U_i$ in $\tilde{H}$ defined by
 \eqspl{U_i}{U_i=\{ u_j=1, \ \forall j< i;\ v_j=1, \ \forall j\geq i\},} and these
$U_i, i=1,..., m$ cover a neighborhood of the special fibre of
$\tilde{H}.$ Now the generators $F_i$ admit the following
relations:
\beq u_{i-1}F_j=u_jx^{i-1-j}F_{i-1},\ 0\leq j<i-1;\ v_iF_j=v_jy^{j-i}F_i,\
m\geq j>i \eeq where we set $u_i=v_i=1$ for $i=0,m.$ Hence $\I$ is
generated on $U_i$ by $F_{i-1}, F_i$ and assertions (iii), (iv)
follow directly from Theorems 1,2 and 3 of \cite{Hilb} .\par
As for (v), it follows immediately from the definition of the
$F_i$,  plus
the fact just noted that, over $U_i,$ the ideal $\I$ is generated by
$F_{i-1}, F_i$, and that on $U_i$, we can set
$u_{i-1}=v_{i}=1.$ Finally (vi) is contained in \cite{Hilb}, Thm. 2.

\end{proof}
At this point it's worth noting the following consequences of Theorem \ref{Hilb-local}, (i). First, recall that a deformation $X/B$ of a nodal curve $X_0$ is said to be \emph{locally versal}
(or \emph{locally versal at the nodes}) if the natural map of $B$ to
the product of local deformation spaces is smooth.
\begin{cor}\label{hilb-is-smooth-for-versal-cor}Let $X/B$ be a family of nodal or smooth curves.
\par
(i) $X\sbr m._B/B$ is a normal crossings morphism, i.e. fibres have normal crossings.
\par (ii) If $X/B$ is locally versal at the nodes, then $X\sbr m._B$
 and the universal subscheme over $X\sbr m._B$ are smooth.
\par (iii) If $X$ is irreducible then so is $X\sbr m._B$\end{cor}
\begin{remark*}
In (ii), the smoothness claimed is of course in the absolute sense, i.e. over $\C$, not over $B$.
\end{remark*}
\begin{proof}We first prove (ii) as (i) is similar and simpler.
Working near a fibre $X_0$,
there is a standard coordinate neighborhood $U_i$ of each node $p_i, i=1,...,k$,
which is a pullback of $V/T: xy=t$, and such that the product map
$B\to  T^k$ is smooth. Then ${\prod}_B(U_i)\sbr m_i._B$ is smooth over
${\prod}_\C V\sbr m_i._T$, and the latter is smooth. Therefore
${\prod}_B(U_i)\sbr m_i._B$ is smooth, hence so is $X\sbr m._B$.
\par (iii) It follows from the local models that the every fibre component of $X/B$
is $m$ dimensional and dominates a fibre component of $X\spr m._B$.
Since $X\spr m._B$ is irreducible, so is $X\sbr m._B$.
\end{proof}

 In light of Theorem \ref{Hilb-local}, we
identify a neighborhood $H_m$ of the special fibre in $\tilde{H}$
with a neighborhood of the punctual Hilbert scheme (i.e. $ c_m\inv(mp)$)
in $X\sbr m._B$, and note that the projection
$H_m\to \A^{2m}\times B$ coincides generically, hence everywhere,
with $\sigma\circ  c_m$. Hence $H_m$ may be viewed as the
subscheme of $U\spr m._B\times_B\tilde{C}$ defined by the
equations
\begin{eqnarray}\label{sigma-uv-rel}\nonumber \sigma_m^xu_1=tv_1, \\
\sigma_{m-1}^xu_1=\sigma_{1}^yv_{1},...,
\sigma^x_{1}u_{m-1}=\sigma^y_{m-1}v_{m-1},\\
\nonumber
 tu_{m-1}=\sigma^y_mv_{m-1}. \end{eqnarray}
Alternatively, in terms of the $Z$ coordinates, $H_m$ may be defined as
the subscheme of
$U\spr m._B\times \P^{m-1}\times B$ defined by the relations
\eqref{Z-rel-quadratic}, which define $\tilde{C}$, together with
\eqspl{sigma-Z-rel}{\sigma_{i}^yZ_i=
\sigma_{m-i}^xZ_{i+1},\
 i=1,...,m-1.}

\newsection{Reverse engineering and proof of Blowup Theorem}
\label{reverse-engineering-sec}

Reverse-engineering an ideal means finding generators with given
syzygies.
 Our task now is effectively
to reverse-engineer an ideal (discriminant ideal) in the $\sigma$'s whose
syzygies, for suitable generators,
are given by \eqref{sigma-Z-rel} and \eqref{Z-rel-quadratic}.
This will be achieved by passing to the ordered version of Hilb,
i.e. $X\sbr m._B\times_{X\spr m._B}X^m_B$.
The sought-for generators will be given by certain 'mixed Van der
Monde'
determinants. The proof of the Blowup Theorem is then concluded,
essentially by showing explicitly that, locally over Hilb,
all the generators are multiples of one of them.
\subsection{Order}
Let
$OH_m=H_m\times_{U\spr m._B}U^m_B$, so we have a cartesian
diagram \beq \begin{matrix} OH_m&\stackrel
{\varpi_m}{\longrightarrow}&
H_m\\ o c_m \downarrow &\square&\downarrow  c_m\\
X^m_B&\stackrel{ \omega_m}{\longrightarrow}& X\spr
m._B\end{matrix}
 \eeq and its global analogue
\beq\begin{matrix} X^{\lceil m\rceil}_B&\stackrel
{\varpi_m}{\longrightarrow}&
X\sbr m._B\\ o c_m \downarrow &\square&\downarrow  c_m\\
X^m_B&\stackrel{ \omega_m}{\longrightarrow}& X\spr
m._B\end{matrix}\eeq
Here the horizontal maps are all $\frak S_m$-quotients, hence flat.
Note that $ X\spr m._B$ is normal and Cohen-Macaulay: this follows from
the fact that it is a quotient by
$\frak S_m$ of $X^m_B$, which is a locally complete intersection
with singular locus of codimension $\geq 2$ (in fact, $>2$, since
$X$ is smooth). Alternatively, normality of $ X\spr m._B$ follows
from the fact that $H_m$ is smooth and the fibres of $ c_m:H_m\to
 X\spr m._B$ are connected and reduced (being products of connected chains of
rational curves), using the following general fact: if $f:A\to B$
is a proper surjective morphism with connected reduced fibres between
integral algebraic schemes over an algebraically closed field 
and $A$ is normal, then so is $B$ [proof: For
any closed point $b\in B$, the inclusion $\O_{B, b}/\m_{B,b}\to H^0(\O_{f\inv(b)})=H^0(\O_A/\m_{B,b}\O_A)$
 is an isomorphism because $f\inv(b)$ is reduced and connected.
By an easy composition series argument, the analogous statement holds with $\m_{B,b}$
replaced by $\m_{B,b}^n$ for any $n\geq 1$. Consequently,
by the Formal Function Theorem (\cite{hart},\S 3.11), we have
$f_*(\O_A)=\O_B$. Then since the local rings of $A$ are integrally closed, 
the same is true of $\O_A(U), U\subset A$ open, hence
also for the local rings of $\O_B$].\par 
Now a few remarks are in order.
\begin{itemize}\item
The map $\omega_m$ is simply ramified
generically over $D^m$ and we have \beq \omega_m^*(D^m)=2OD^m
\eeq where
\beq OD^m=\sum\limits_{i<j}D^m_{i,j} \eeq where
$D^m_{i,j}=p_{i,j}\inv(OD^2)$ is the locus of points whose $i$th and
$j$th components coincide.\par
\item
$\Bl_{OD^m}X^m_B=\Bl_{2OD^m}X^m_B$ as blowups,  because
blowing up an ideal and its powers are the same (see \cite{hart}, Ex.
II.7.11.a).
\item We have
\[(\Bl_{D^m}X\spr m._B)\times_{X\spr m._B}X^m_B=
\Bl_{\omega_m^*(D^m)}X^m_B=\Bl_{2OD^m}X^m_B.\]
So if the natural map $X\scl m._B\to\Bl_{2OD^m}X^m_B$ is an
isomorphism,
then (obviously) so is the $\frak S_m$-equivariant map
\[f: X\scl m._B\to (\Bl_{D^m}X\spr m._B)\times_{X\spr m._B}X^m_B,\]
which is just the pullback of the natural map
\[c_m': X\sbr m._B\to\Bl_{D^m}X\spr m._B\]
by the finite flat surjective map $\varpi_m$, therefore so is
$c_m'$ itself (which is the $\frak S_m$-quotient of $f$).
\item Therefore finally we are reduced to showing that $oc_m$ is
    equivalent to the blowing up of $OD^m$.
\end{itemize}

\par

\par
The advantage of working with $OD^m$ rather than its unordered
analogue is that at least some of its equations are easy to write
down: let \beq v^m_x=\prod_{1\leq i<j\leq m}(x_i-x_j), \eeq and likewise
for $v^m_y.$ As is well known, $v^m_x$ is the determinant of the Van
der Monde matrix
\beq V^m_x=\left [\begin{matrix}
1&\ldots&1\\x_1&\ldots&x_m\\\vdots&&\vdots\\
x_1^{m-1}&\ldots&x_m^{m-1}\end{matrix} \right ]. \eeq
 Also
set
\beq \tilde{U}_i=\varpi_m\inv(U_i), \eeq where $U_i$ is as in \eqref{U_i},
being a neighborhood of $q_i$ on $H_m$. Then in $U_1$, the
universal ideal $\I$ is defined by
\beq F_0, \ \ F_1=y+(\text{function of }x) \eeq and consequently the
length-$m$ scheme corresponding to $\I$ maps isomorphically to its
projection to the $x$-axis. Therefore over
$\tilde{U}_1=\varpi_m\inv(U_1),$ where $F_0$ splits as $\prod
(x-x_i),$ the equation of $OD^m$ is simply given by \beq G_1=v^m_x.
\eeq
Similarly, the equation of $OD^m$ in $\tilde{U}_m$ is given by
\beq G_m=v^m_y. \eeq Now let \beq \Xi:OH_m\to\P^{m-1} \eeq be the
morphism
corresponding to (the pullback of) $[Z_1,...,Z_m]$ (cf. \eqref{Z_i}; note that this is by definition
essentially a product projection,
hence a morphism); set \eqspl{L-def}
{L=\Xi^*(\O(1)).} Note
that $\tilde{U}_i$ coincides with the open set where $Z_i\neq 0$,
so $Z_i$ generates $L$ over $\tilde{U}_i.$ Let
\beq O\Gamma\spr m.=o c_m\inv(OD^m). \eeq
We shall see below that this is a Cartier divisor, in fact we
shall construct an isomorphism
\beql{gamma}{ \gamma:\O(-O\Gamma\spr m.)\to L.}
This isomorphism
will easily imply Theorem 1. To construct $\gamma,$ it suffices
to specify it on each $\tilde{U}_i.$
\subsection{ Mixed Van der
Mondes}  A clue as to how the
latter might be
done comes from the relations \eqref{sigma-Z-rel}. Thus, set
\beql{Gdef1}{ G_i=\frac{(\sigma_m^y)^{i-1}}{t^{(i-1)(m-i/2)}}v^m_x=
\frac{(\sigma_m^y)^{i-1}}{t^{(i-1)(m-i/2)}}G_1,\ \ i=2,...,m.} Thus,
\beql{Gdef2}{ G_2=\frac{\sigma^y_m}{t^{m-1}}G_1,
G_3=\frac{\sigma^y_m}{t^{m-2}}G_2,...,
G_{i+1}=\frac{\sigma^y_m}{t^{m-i}}G_i, i=1,...,m-1.}
In light of \eqref{sigma-x-y-rel},\eqref{sigma-x-y-rel2}, we deduce
\eqspl{sigma-G-rel}{
\sigma_{m-i}^xG_{i+1}=\sigma_i^yG_i.
}Comparing this with \eqref{sigma-Z-rel}
certainly suggests solving our reverse-engineering problem by
 assigning $Z_i$ to $G_i$, which is what we will do eventually.
\begin{rem}
Clearly $v_x^m$, hence all the $G_i$, are invariant under the
alternating group $\frak A_m$, hence are well-defined on the
'Orientation product', i.e. the quotient of $X^m_B$ by the action
of $\frak A_m$, which coincides with the double cover of $X\spr m._B$
branched on $D^m$.
\end{rem}
Now an elementary calculation shows that if we denote by
$V^m_i$ the 'mixed Van der Monde' matrix
\beq V^m_i=\left[\begin{matrix}
1&\ldots&1\\x_1&\ldots&x_m\\\vdots&&\vdots\\
x_1^{m-i}&\ldots&x_m^{m-i}\\y_1&\ldots&y_m\\\vdots&&\vdots\\
y_1^{i-1}&\ldots&y_m^{i-1}\end{matrix}\right]\eeq then we
have
\beql{vandermonde}{ G_i=\pm\det (V^m_i), i=1,...,m.}
\begin{minipage}[0pt,50pt]{450pt}
$\ulcorner$ Indeed for $i=1$ this is standard; for general $i$, it suffices
to prove the analogue of \refp{Gdef2}. for the mixed Van der Monde
determinants. For this, it suffices to multiply each $j$th column of $V^m_i$
by $y_j$,
and factor a $t=x_jy_j$ out of each of rows  $2,...,m-i+1$,
 which yields
\beql{vandermonde1}{\sigma^y_m\det(V^m_i)=(-1)^{m-i}t^{m-i}V^m_{i+1}.}
$\llcorner$\end{minipage}\par
    From \refp{vandermonde}. it follows, e.g., that $G_m$ as given in
    \refp{Gdef1}. coincides with
$\pm v^m_y.$
\subsection{Conclusion of proof} The following
result is key for the Blowup Theorem. \begin{lem}
\label{Gi-generates-on-Ui} $G_i$ generates $\O(-O\Gamma\spr m.)$
over
$\tilde{U}_i.$ In particular, $O\Gamma\spr m.$ is Cartier.
\end{lem}\begin{proof}[Proof of Lemma] This is clearly true where $t\neq
0$ and it remains
to check it along the special fibre $OH_{m,0}$ of $OH_m$ over $B$.
Note that $OH_{m,0}$ is a sum of components of the form
\beql{Theta} {\Theta_I=\text{Zeros}(x_i,i\not\in I, y_i, i\in I),
 I\subseteq\{1,...,m\}, } none of which is contained in the
 singular locus of $OH_m.$
  Set
 \beq \Theta_i=\bigcup\limits_{|I|=i}\Theta_I. \eeq Note that
 \beq \Tilde{C}_i\times 0\subset \Theta_i, i=1,...,m-1 \eeq and therefore
 \beq \tilde{U}_i\cap\Theta_j=\emptyset, j\neq i-1, i. \eeq
 Note that $y_i$ vanishes to order 1 (resp. 0) on $\Theta_I$
 whenever $i\in I$ (resp. $i\not\in I$). Similarly, $x_i-x_j$
 vanishes to order 1 (resp. 0) on $\Theta_I$ whenever both $i,j\in
 I^c$ (resp. not both $i,j\in I^c$). From this, an elementary
 calculation shows that the vanishing order of
 $G_j$ on every component $\Theta$ of
 $\Theta_k$ is \beql{ordG}{ \ord_{\Theta}(G_j)=(k-j)^2+(k-j).}
 We may unambiguously denote this number by $\ord_{\Theta_k}(G_j)$.
 Since this order is nonnegative for all $k,j,$
 it follows firstly that
 the rational function $G_j$ has no poles, hence is in fact regular
 on $X^m_B$ near
 $mp$ (recall that $X^m_B$ is normal); of course,
 regularity of $G_j$ is also immediate from \refp{vandermonde}..
  Secondly, since this order is zero for $k=j, j-1$, and
 $\Theta_j, \Theta_{j-1}$ contain all the components of $OH_{m,0}$
  meeting
 $\tilde{U}_j$, it follows that in $\tilde{U}_j,$ $G_j$ has no
 zeros besides $O\Gamma\spr m.\cap\tilde{U}_j,$ so $G_j$ is a
 generator
 of $\O(-O\Gamma\spr m.)$ over $\tilde{U}_j.$ QED Lemma.
  \end{proof}
  The Lemma yields a set of generators for the ideal of $OD^m$:
  \begin{cor}[of Lemma]\label{Gs-generate-bigdiag}
  The ideal of $OD^m$ is generated, locally
 near $(p,...,p)$, by $G_1,...,G_m.$  \end{cor}\begin{proof}
  If $Q$ denotes the cokernel of the map
  $\bigoplus\limits_m\O_{X^m_B}\to \O_{X^m_B}(-{OD^m})$ given by
  $G_1,...,G_m$,
  then $oc_m^*(Q)=0$ by the Lemma,
   hence $Q=0$, so the $G$'s generate $\O_{X^m}(-{OD^m})$.
  \end{proof}
  Now we can construct the desired isomorphism $\gamma$ as in
  \eqref{gamma}, as
   follows. Since $Z_j$ is a
 generator of $L$ over $\tilde{U}_j,$ we can define our
 isomorphism $\gamma$ over $\tilde{U}_j$ simply by specifying that
 \beq \gamma (G_j)=Z_j\  {\text{on}}\ \tilde{U}_j. \eeq Now to check that
 these maps are compatible, it suffices to check that
 \beq G_j/G_k=Z_j/Z_k \eeq as rational functions (in fact, units over
 $\tilde{U}_j\cap\tilde{U}_k$). But the ratios $Z_j/Z_k$ are
 determined by the relations \refp{sigma-Z-rel}., while $G_j/G_k$
 can be computed
 from \refp{sigma-G-rel}., and it is trivial to check that these agree. \par
 Now we can easily complete
 the proof of Theorem 1. The existence of $\gamma$, together with the
 universal property of blowing up, yields a morphism
 \[ Bc_m:OH_m\to B_{OD^m}X^m_B\]
 which is clearly proper and birational, hence surjective.
 On the other hand, the fact that the $G$'s generate the ideal
 of $OD^m$, and correspond to the $Z$ coordinates on $OH_m\subset
 X^m_B\times \P^{m-1}$, implies that $Bc_m$
 is a closed immersion. Therefore $Bc_m$ is an isomorphism.
 \qed
 \begin{rem}
 It follows from the foregoing proof that the cycle map
 is a \emph{nonlinear} blowup, i.e. that
 the inclusion $\Proj(\bigoplus I_{D^m}^n)\subset\P(I_{D^m})$
 is proper.
 \end{rem}
 \subsection{Complements and consequences}
 These concern the standard family $V/T$ given by $xy=t$:
 \begin{cor}\label{sym-eq} For the standard family, the image of the
 relative symmetric product $V\spr m._T$ under the elementary
 symmetric functions embedding $\sigma$ (cf. Lemma \ref{sigma-emb}) is
 schematically defined by the equations
 (\ref{sigma-x-y-rel}-\ref{sigma-x-y-rel2}).\end{cor}
 \begin{proof}
 We have a diagram locally
 \beql{}{\begin{matrix}
 H^m&\subset& \P^{m-1}\times\A^{2m}\times T\\
 \downarrow&&\downarrow\\
 V\spr m._T&\stackrel{\sigma}{\hookrightarrow}&\A^{2m}\times T.
 \end{matrix}
}
We have seen that the image of the top inclusion is defined by the
equations \eqref{Z-rel-quadratic}, \eqref{sigma-Z-rel}.
 The equations of the schematic image of $\sigma$ are obtained
 by eliminating the $Z$ coordinates from the latter equations,
 and this clearly yields the equations as claimed.
 \end{proof}
 Now as one byproduct of the proof of Theorem \ref{blowup}, we obtained
 generators of the ideal of the ordered half-discriminant $OD^m$.
As a further consequence, we can determine the ideal of the
 discriminant locus $D^m$ in the symmetric product $X\spr m._B$ itself:
 let $\delta_m^x$
 denote the discriminant of $F_0$, which, as is
 well known \cite{L}, is a polynomial in the $\sigma_i^x$
 such that \beql{}{ \delta_m^x=G_1^2.} Set
 \beql{eta}{\eta_{i,j}=\frac{(\sigma_m^y)^{i+j-2}}
 {t^{(i-1)(m-i/2)+(j-1)(m-j/2)}}
 \delta^m_x.} It is easy to see that this is a polynomial in the $\sigma^x_.$
 and the $\sigma^y_.$, such that $\eta_{i,j}=G_iG_j.$
 \begin{cor} For the standard family $xy=t$, the ideal of $D^m$ is generated, locally
 near $mp$, by $\eta_{i,j}, i,j =1,...,m.$\end{cor}\begin{proof} This
 follows from the fact that $\varpi_m$ is flat (being a $\frak S_m$-
 quotient) and that
 \beq \varpi_m^*(\eta_{i,j})=G_iG_j, i,j=1,...,m \eeq generate the ideal of
 $2OD^m=\varpi_m^*(D^m).$\end{proof}
 Because any family of nodal curves $X/B$ is locally isomorphic to a pullback
 of the standard family $V/T$ , it follows that analogues of the previous two
 corollaries hold for $X\spr m._B$ over a neighbourhood of a point
 $m\theta$, where $\theta$ is a relative node.

 \newsection{Discriminant polarization}
 We now return to the case of a general family $X/B$ of nodal-or-smooth
 curves.
 We study some natural sheaves, including the discriminant polarization
 , on
 the Hilbert schemes $X\sbr m._B$.
   \par
 Note that the ideal of the Cartier divisor $ c_m^*(D^m)$ on $X\sbr m._B$,
 that is,
 $\O_{X\sbr m._B}(- c_m^*(D^m))$,
 is isomorphic in terms of our local model $\tilde{H}$ to
 $\O(2)$ (i.e. the pullback of $\O(2)$ from $\P^{m-1}$).
 This suggests that $\O(- c_m^*(D^m))$
  is divisible by 2 as line bundle on $X\sbr m._B$.
  This is indeed so, and is subsumed in the definition of \emph{discriminant
polarization} which follows, together with that of \emph{tautological sheaf}. 
Consider the tautological subscheme
\[D^{m,1}\subset X\sbr m._B\times_BX\]
with maps $p_X:D^{m,1}\to X, p_{X\sbr m._B}:D^{m,1}\to X\sbr m._B$.
 \begin{defn}(i)
For any sheaf $A$ on $X$, the associated tautological sheaf is
defined by
\[\Lambda_m(A)=p_{X\sbr m._B*}(p_{X}^*(A))\]
(ii)  The discriminant polarization on $X\sbr m._B$ is defined by
\eqsp{
\O_{X\sbr m._B}(1)=\O(-\Gamma\spr m.):=\det(\Lambda_m(\O_X))
}
\end{defn}
Note that if $A$ is locally free, then by flatness so is $\Lambda_m(A)$.
These bundles are obviously compatible with base-change.
Moreover, note that the trace pairing \[\Lambda_m(\O_X)\otimes
\Lambda_m(\O_X)\to \O_{X\sbr m._B}\] yields
a generically injective map \mbox{$\Lambda_m(\O_X)\to \Lambda_m(\O_X)^*$}
which drops rank precisely
on the discriminant $c_m^*(D^m)$, therefore $2\Gamma\spr m.\sim_{\mathrm{lin}}
c_m^*(D^m)$.
\par 
 We will also use the notation
 \beq\O(\Gamma\spr m.)=\O_{X\sbr m._B}(-1).\eeq
Note 
 that $\Gamma\spr m.$ is defined as an effective Weil divisor, and as a line bundle, but not necessarily 
 as an effective Cartier divisor,
 though $2\Gamma\spr
 m.$ and $\Gamma\scl m.$ are effective (the latter because the
 symmetrization map $X^m_B\to X\spr m._B$ is generically ramified with
 multiplicity 2 along $D^m$).
 In fact, $\Gamma\spr m.$ is essentially never effective Cartier, as
 Remark \ref{divisible} below shows.
 Nonetheless, $-\Gamma\spr m.$
  is relatively ample on the Hilbert
 scheme $X\sbr m._B$ over the symmetric product $X\spr m._B$,
hence the name {discriminant polarization.}

Further light on the discriminant is shed
by the notion of \emph{norm}:
\begin{defn}
For a line bundle $A$ on $X$, its $m$-th norm on $X\sbr m._B$ is defined
by
\[[m]_*(A)=\det(\Lambda_m(A))\otimes\O(\Gamma\spr m.)\]
\end{defn}
If $A=\O(Y)$ for an effective divisor $Y$, the exact
sequence
\[0\to \Lambda_m(A^*)\to \Lambda_m(\O_X)\to \Lambda_m(\O_Y)\to 0\]
shows that in this case $[m]_*(A)=-[m]_*(A^*)=\det(\Lambda_m(\O_Y))$ is
an effective
divisor supported on the locus of schemes whose support
meets that of $Y$.

 \begin{rem}\label{divisible}
 Let $X$ be a smooth curve of genus $g\geq 2$ and fix $m\geq 2$.
 Then the discriminant $D\subset X\spr m.$ is not algebraically
 equivalent
 to $\sum a_iA_i$ where each $a_i>0$, $\sum a_i\geq 2$ and the $A_i$
 are effective and nontrivial; thus, $D$ is neither splittable
 nor divisible as effective divisor up to algebraic equivalence.\par
\noindent\emph{Proof.}\quad Else, it follows that  $D$, being a prime
divisor,
 meets each $A_i$ properly, hence
 $\O_D(A_i)$ is effective, therefore $\O_D(D)$ is effective up to
 algebraic equivalence on $D$. Letting $f:X\times X\spr m-2.\to D$
denote the obvious (normalization) map, $f(x,z)=2x+z$, it follows
that $f^*(D)$ is effective. But
 \[f^*(D).(X\times{\mathrm{pt}})=-\deg(\omega_X)=2-2g<0,\] which
contradicts effectivity.\par
For $g\leq 1$, $D$ is effectively divisible by 2, at least for a
single curve. For $g=1$,
$X$ is an elliptic curve with group law $*$ and
$D$ is algebraically equivalent to $2D_a, a\in X$, where
\[D_a=\{x+x*a+\sum\limits_{i=1}^{m-2}x_i\}\]
The algebraic equivalence becomes linear when $a$ has order $2$
in the group.

 \end{rem}
  \newsection{Flags}\label{flags-sec}
 See \cite{sernesi} for Flag- Hilbert schemes in general. Flag-Hilbert schemes for points on nodal curves
  were studied in \cite{geonodal, Hilb}. In \cite{geonodal}, a construction is
  given for the full-flag
 Hilbert scheme via an explicit blowup procedure, different in
 flavor from the above discriminant blowup.
 In \cite{Hilb}, a model analogous to $H_m$ was constructed for the
 relative Hilbert scheme $X\sbr m,m+1._B$ of $(m,\nobreak m+\nobreak
 1)$-flags,
 i.e. pairs of ideals $(z_1\supset z_2)$ of respective lengths
 $(m,m+1)$.
 Here we try to reconcile the two viewpoints by
 showing that the full-flag Hilbert scheme can also be represented
 as a blowup of a discriminant-like (viz. incidence) variety, in analogy with the case of the ordinary Hilbert scheme.\par
 Consider the flag-Hilbert scheme, which fits in a diagram
 \eqspl{}{\begin{matrix}
 &&X\sbr m,m+1._B&\subset& X\sbr m._B\times X\sbr m+1._B\\
 &\ \ \ \ \ \ \ \ \ p_{[m]}\swarrow&\downarrow p_{[m+1]}&&\\
 &X\sbr m._B\ \ \ \ &X\sbr m+1._B
 \end{matrix}
 } Via this, $X\sbr m,m+1._B$ is endowed with divisors
   denoted $\Gamma\spr m., \Gamma\spr m+1.$, which are
   pullbacks of analogous divisors on $X\sbr m._B, X\sbr m+1._B$
   respecively.
   There is a natural morphism (where $X$ is identified with the set of
   colength-1 ideals)
 \eqsp{X\sbr m,m+1._B\to X\\
 (z_1\supset z_2)\mapsto \Ann(z_1/z_2)
 } whence a map
 \eqspl{}{
 c_{m,1}:X\sbr m,m+1._B\to X\sbr m._B\times_BX
 }
 \begin{thm}\label{blowup-flag-thm}
 $c_{m,1}$ is the blowing-up of the incidence variety
 $D\spr m,1.=\{(z, x): x\in z\}$
 \end{thm}
 \begin{proof}
 Let \[b:Y\to X\sbr m._B\times_BX\] be the blowing up of $D\spr m,1.$ and
 $\Gamma\spr m,1.$ the
 exceptional (Cartier) divisor, i.e. the inverse image of $D\spr m,1.$.
 Because $c_{m,1}\inv(D\spr m,1.)=\Gamma\spr m+1.-\Gamma\spr m.$
 is
 Cartier, it follows from the universal property of blowing up that we get a
 diagram
 \eqsp{\begin{matrix}
 X\sbr m,m+1._B&&\stackrel{c'}{\to}&& Y\\
 c_{m,1}&\searrow&&\swarrow b&\\
 &&X\sbr m._B\times_BX&&
 \end{matrix}
 }
 On the other hand, there is an obvious map
 \[Y\to X\spr m+1._B\] and the pullback of $D\spr m+1.$ is just
 $\Gamma\spr m.+\Gamma\spr m,1.$, hence Cartier. So by the Blowup
 Theorem we
 get a map $Y\to X\sbr m+1._B$. Together with the projection
 $Y\to X\sbr m._B$, this gives a map
 $Y\to X\sbr m._B\times _BX\sbr m+1._B$ whose image is clearly
 contained in  $X\sbr m,m+1._B$ , whence
 a map \[d: Y\to X\sbr m,m+1._B,\] which together with $c'$ fits in  a
 diagram
 \eqsp{\begin{matrix}
 X\sbr m,m+1._B&&\stackrel[d]{c'}{\rightleftarrows}&& Y\\
 c_{m,1}&\searrow&&\swarrow b&\\
 &&X\sbr m._B\times_BX&&
 \end{matrix}
 }
 
 As both vertical maps are birational, $c',d$ are
 mutually inverse isomorphisms.
 \end{proof}
 As a consequence, we obtain recursively a presentation of the full-flag
 Hilbert scheme as a blowup of incidence varieties. This slightly generalizes a result proven
 in (\cite{geonodal}, Thm. 2.1) by more explicit means.
 \begin{cor}\label{blowup-flag-cor}
 Denote by $W^{m.}(X/B)$ the flag Hilbert scheme parametrizing
 flags of subschemes of fibres $(z_{m_1}<z_{m_2}...<z_{m_k})$
 of respective lengths $m_1<m_2<...<m_k$. Then
 $W^{m., m_k+1}(X/B)$
 is the blowup of $W^{m.}(X/B)\times_BX$ in the incidence variety
 \[D\spr{m.},1.=\{(z.,x):x\in z_{m_k}\}.\]
 \end{cor}
 \begin{rem}
 We don't know if the analogues of Theorem \ref{blowup-flag-thm} or Corollary \ref{blowup-flag-cor}
 hold for arbitrary flags, e.g. of type $[m, m+2]$. Those Hilbert schemes seem to be worse behaved:
 inter alia,  the fibres of the cycle map on $X^{[m,m+2]}$ can have dimension 2 if $m>1$. For instance,
 a generic length-2 subscheme of a node is contained in just two length-3  subschemes, but in an
 entire 1-paramater family of length-4 subschemes.
 \end{rem}

\part{Node scrolls}
 \newsection{Study of $H_m$ }\label{study of H_m}
We continue our study of the cycle map over a neighborhood of
a maximally singular cycle $m\theta$ with $\theta$ a fibre node, using
the model $H_m$. The results will be applied in the Node Scroll
theorem.
Having previously determined the structure of $c_m$ along its 'most
special' fibre $c_m\inv(m\theta)$ (which corresponds in the model
$H_m$
 to the fibre over the origin $0_{\A^{2m}}$), our purpose in this
section is to
determine its structure along nearby fibres and their variation. Thus we will
assume for the rest of this
section, unless otherwise stated, that we are in the
local situation where $B$ is a
smooth curve, with local coordinate $t$, and the family
$U/B$ is the standard degeneration $xy=t$. Our purpose is to prove
the following result, which serves as the foundation for our study of node scrolls.
The notation will be explained below; suffice it to recall here
that on a node with equation $xy=0$, an ideal of type $C^n_j$ (resp. $Q^n_j$)
is generated by $x^{n-j}+ty^j, t\neq 0$ (resp. $x^{n-j+1}$ and $ y^j$).
\begin{lem}\label{nodescroll-lem}
For each $1\leq j\leq n-1$, there exists a $\P^1$-bundle $F^{m,n}_j$
over $(U^\theta)\spr m-n.$, together with a pair of disjoint
sections $Q^{m,n}_j, Q^{m,n}_{j+1}$ and a map
\[p_{j,[m]}:F^{m,n}_j\to H_m,\]
such that
\begin{enumerate}\item
the image of $p_{j,[m]}$ coincides with the closure of the locus
of schemes having length $n$ and type $C^n_j$ at $\theta$;
\item the combined image of
\[\coprod \limits_{j=1}^{n-1} F^{m,n}_j\to H_m\]
coincides with the locus of schemes of length at least $n$ at $\theta$
\item  the image  $p_{j,[m]}(Q^{m,n}_\bullet), \bullet=j,j+1$, coincides
    with
the closure of
the locus of schemes having length $n$ and type $Q^n_\bullet$ at
$\theta$.
\end{enumerate}

\end{lem}
  \subsection{Nearby fibres} Let $U', U"$ denote the $x,y$ axes,
  respectively
in $U_0=X_0\cap U$, with their respective origins $\theta', \theta"$
mapping to $\theta\in U$. Set $U^\theta=U'\coprod U"$, the normalization
of $U_0$. If the
special fibre $X_0$ is reducible, then $U', U"$ globalize to
(i.e. are open subsets of) the two
components of the normalization. If $X_0$ is irreducible, then
both $U'$ and $U"$ globalize to the normalization.
 For any pair of natural numbers $(a,b), 0<a+b< m$,
set $$U^{(a,b)}=U^{'(a)}\times U^{"(b)}$$
(which globalizes to a component --the unique one,
if $X_0$ is irreducible-- of the normalization of $X_0^{a+b}$).
 Then we have a
natural map $$U^{(a,b)}\to(U_0)\spr m._B\subset
(U)\spr m._B$$ given by
\beq (\sum m_ix_i,\sum n_jy_j)\mapsto \sum
m_i(x_i,0)+\sum n_j(0,y_j)+(m-a-b)\theta. \eeq  This map
is clearly birational to its image, which we denote by $\bar U^{(a,b)}.$ Thus
$U^{(a,b)}$  coincides with the normalization of $\bar U^{(a,b)}$.
 It is clear that $\bar U^{(a,b)}$ is
defined by the equations \beq \sigma^x_m=...=\sigma^x_{a+1}=0,
\sigma^y_m=...=\sigma^y_{b+1}=0. \eeq
 A  point
$$c\in \bar U^{(a,b)}-(\bar U^{(a+1,b)}\cup \bar U^{(a,b+1)}),$$ i.e. a
cycle in which $(0,0)$ appears
with multiplicity exactly $n=m-a-b$, is said to be of \emph{type $(a,b)$}.
Type yields a natural stratification of the symmetric product $U\spr
m._0$.
Now let $\bar H^{(a,b)}$ be the closure of the locus of schemes whose
cycle is of type $(a,b)$, i.e.
\beql{Hab}{\bar H^{(a,b)}={\mathrm{closure}}(c_m\inv(\bar U^{(a,b)}-(\bar
U^{(a+1,b)}\cup \bar U^{(a,b+1)})))\subset H_m}
Also let
\eqspl{}{
H\spr a,b.=\bar H\spr a,b.\times_{\bar U\spr a,b.}U\spr a,b. .
}
Clearly the restriction of $c_m$ on $\bar H^{(a,b)}$ factors through a
map
\begin{eqnarray*}\tilde c_m: \bar H^{(a,b)}\to \bar U^{(a,b)},\\
\tilde c_m=((\sigma^x_1, ...,\sigma^x_a), (\sigma^y_1,
...,\sigma^y_b))\end{eqnarray*}
Approaching the 'origin cycle' $m(0,0)$ through
cycles of type $(a,b)$,
on $\bar U^{(a,b)}$, means that $a$ (resp. $b$) points are approaching
the origin $\theta'$ (resp. $\theta"$) along the $x$ (resp. $y$)-axis.  For a
general cycle $c$ of type $(a,b)$, we have, for all $j\leq b$, that
$\sigma^y_j\neq 0, \sigma_{m-j}^x=0$,  hence
by the equations \eqref{H-equations}
(setting each $a_i=\sigma^x_{m-i},
d_i=\sigma^y_{m-i}$), we conclude $v_j=0$; thus
\beql{v=0}{ v_1=...=v_{b}=0; } similarly, for all $j\leq a$,
we have $\sigma_{m-j}^y=0, \sigma_j^x\neq 0$ ($c$ being general),
hence
again by the equations \eqref{H-equations} , we conclude $u_{m-j}=0$;
thus
\beql{u=0}{  u_{m-1}=...=u_{m-a}=0. }
Consequently, the fibre of $ c_m $ over this point is schematically
\beql{cm fibre}{c_m\inv(c)=\tilde c_m\inv(c)\simeq
\bigcup\limits_{i=b+1}^{m-a-1} C^m_i, } provided
$a+b\leq m-2$ (where the $C^m_i$ are the
 components of the punctual Hilbert scheme, as in the basic construction
of the model $H_m$, see Theorem \ref{Hilb-local}). If $a+b=m-1$, the fibre
is the unique point given by
\beq v_1=...=v_b=u_{b+1}=...=u_{m-1}=0  \eeq (as a subscheme
of $X/B$, this point is the
one denoted
$Q^m_{b+1}$ in \cite{Hilb}, and has ideal $(x^{m-b}, y^{b+1})$). As $c$
approaches the 'origin' $(m\theta)$ in $\overline{U}^{(a,b)}$,
 or for that matter any point $c'$,
 the equations \refp{v=0}.,\refp{u=0}. persist, so we conclude
\beql{tilde cm fibre} {\tilde
c_m\inv(c')=\begin{cases}\bigcup\limits_{i=b+1}^{m-a-1} C^m_i, a+b\leq
m-2,\\
Q^m_{b+1}, a+b=m-1.
\end{cases}
}
[Informally, this is a priori plausible: because schemes in $C_i^m$
represent $i$ points coalesced through the $y$-axis and $m-i$
points coalesced through the $x$-axis. Then moving 'out'
 to $c$ represents
generalizing $b<i$ (resp. $a<m-i$) of the $i$
(resp. $m-i$) points over the $y$ (resp. $x$) axis.]\par
Of particular interest naturally is the case where the union
above is a single $\P^1$, in other words when $b=i-1, a=m-i-1=m-b-2$.
In this case \[\bar H\spr m-i-1, i-1.\to \bar U\spr m-i-1, i-1.\]
is just a
$\P^1$-bundle, with fibre $C^m_i$ at the origin. Of course the
same is true with the bars removed (i.e. after pullback over $U\spr
m-i-1,i-1.$). [Informally again, this says $C^m_i$ as a component of the
punctual Hilbert scheme (schemes of length $m$
 concentrated at $\theta$) extends most generically by freeing up $i-1$
 and $m-i-1$ points
respectively over the two axes.]\par

More generally, for any $1\leq j< n\leq m-1, a+b=m-n$,
we have a natural map
\eqsp{\alpha(n-j-1, j-1): U\spr a,b.\to U\spr a+n-j-1, b+j-1.,\\ (.,.)\mapsto
(.+(n-j-1)\theta', .+(j-1)\theta")}
Pulling back over $H\spr a+n-j-1, b+j-1.$, we obtain  $\P^1$-bundles
\eqspl{}{
F^{m,n}_j(a,b)&\to U\spr a,b.\\
F^{m,n}_j=\coprod\limits_{a+b=m-n}F^{m,n}_j(a,b)&\to (U^\theta)\spr
m-n.=\coprod\limits_{a+b=m-n} U\spr a,b. .} We call $F^{m,n}_j$ a
'model node scroll'. It  is a special case of the general
\emph{node scroll}, to be studied further below. Note that
$F^{m,n}_j$ comes equipped with a map $F^{m,n}_j\to H_m$, whose
combined image for $j=1,..,n-1$ by definition is the closure of the locus of
schemes having
length $n$  at the node $\theta$. Note that any subscheme $z$ having
length $n$ locally at $\theta$ sits over a cycle $c$ of type $(a,b),
a+b=m-n$ and therefore occurs in \eqref{cm fibre} for
some $i$, hence also in in
$F^{m,n}_j$ with $j= i-b$. Furthermore, if $z'$ is a subscheme
having length $n'\geq n$ at $\theta$, it occurs on
$F^{m,n}_j(a',b'), a'+b'=m-n'$ for some $j$.
Then choosing $a\geq a', b\geq b'$
with $a+b=m-n$, we can factor $\alpha(n'-j-1, j-1)$ via $U\spr
a,b.$:
\[U\spr a',b'.\to U\spr a,b. \stackrel{\alpha(n-j-1, j-1)}{\to}
U\spr a'+n'-j-1, b'+j-1.\]
to conclude that $z'$ occurs on $F^{m,n}_j(a,b)$ and in particular
on $F^{m,n}_j$. Thus, the image of $ F^{m,n}_j$ in $H_m$ corresponds
to the closure of the locus of schemes which are of length $n$ and
type $C^n_j$ (i.e. local equation $x^{n-j}+\alpha y^j), \alpha\in\C^*$) at the
node $\theta$.
\par Also, referring to \eqref{tilde cm fibre}, we see that
$F^{m,n}_j(a,b)$ and also $F^{m,n}_j$ contain two special, mutually
disjoint cross sections corresponding to $Q^m_{j}, Q^m_{j+1}$, which
come respectively from
\[\bar H\spr m-i, i-1., \bar H\spr m-i-1, i.\subset \bar H\spr m-i-1,
i-1..\] We denote these by $Q^{m,n}_j(a,b),Q^{m,n}_{j+1}(a,b)$
and $Q^{m,n}_j, Q^{m,n}_{j+1}$,
respectively. This notation is slightly imprecise in that there is a
$Q^{m,n}_j$ on both $F^{m,n}_j$ and $F^{m,n}_{j-1}$. But both of
them have the same image in the Hilbert scheme, viz. the closure of
the locus of schemes having length $n$ and type $Q^{m,n}_j$
 (i.e. local equations $(x^{n-j+1}, y^j)$ at $\theta$.
The reason is the same as given above in the case of $F^{m,n}_j$.
This completes the proof of Lemma \ref{nodescroll-lem}.\par
\subsection{Node scrolls: an optional preview}
This subsection is not needed anywhere.
It presents an alternative, more 'qualitative' perspective on
a property of node scrolls that is subsumed in the Node Scroll Theorem
\ref{nodescroll-thm}. This property has to do with the intrinsic,
as opposed to polarized, structure of these scrolls.\par Fixing
$m,n,a,b$ for now, the $F_j=F_j^{m,n}(a,b)$ are components of
special (but typical) cases of what are to be
 called \emph{node scrolls}.
It follows from Lemma \ref{nodescroll-lem} that we can write
$$F_j=\P(L^n_{j}\oplus L^n_{j+1})$$
for certain line bundles $L^n_j, L^n_{j+1}$ on ${U}^{(a,b)}$,
corresponding to the disjoint sections $Q^{m,n}_j, Q^{m,n}_{j+1}$,
where the
difference $L_j^n-L^n_{j+1}$ is uniquely determined (we use additive
notation for the tensor product of line bundles
and quotient convention for projective bundles). The identification of a
natural choice for both these line bundles,
using methods to be developed later in this section,
 will be taken up in the next section and plays an important role in the
 enumerative
geometry of the Hilbert scheme. But the difference
$L_j^n-L^n_{j+1}$, and hence the intrinsic structure of
the node scroll $F_j$, may already be computed now,  as follows.\par
Write
$$Q_j=\P(L_j), Q_{j+1}=\P(L_{j+1})$$ for the two special sections of type
$Q^{m,n}_j, Q^{m,n}_{j+1}$ respectively. Let
$$D_{\theta'}, D_{\theta"}\subset U^{(a,b)}$$ be the divisors
comprised of cycles containing $\theta'$ (resp. $\theta"$).
In the local model, these are given locally by the
respective equations
$$D_{\theta'}=(\sigma^x_a), D_{\theta"}=(\sigma^y_b).$$
\begin{lem}\label{nodescroll0} We have, using the quotient convention for
projective bundles,
\beql{nodescroll-1}{F_j=\P_{U^{(a,b)}}(\O(-D_{\theta'})\oplus
\O(-D_{\theta"})), j=1,...,n-1.}\end{lem}
\begin{proof} Our key tool is a $\C^*$- parametrized family
of sections
'interpolating' between  $Q_j$
and $Q_{j+1}$. Namely,
note that for any $s\in\C^*$, there is a well-defined
section $I_s$ of $F_j$ whose fibre over a general point
$z\in X^{(a,b)}$ is the scheme
$$I_s(z)=(sx^{n-j}+y^j)\amalg \sch(z),$$
where $\sch(z)$ is the unique subscheme of length $a+b$, disjoint from
the nodes, corresponding to $z$, and we are identifying a (principal) ideal
with the corresponding subscheme..\par
\emph{Claim:} The fibre of $I_s$ over a point $z\in D_{\theta'}$
(resp. $z\in D_{\theta"}$) is a scheme of type
$Q^{m,n}_j$, i.e. $(x^{n-j+1}, y^j)$ (resp. $Q^{m,n}_{j+1}$).\par
\emph{Proof of claim}. Indeed set-theoretically the claim is clear from the
fact thar this fibre corresponds to a length-$n$ punctual scheme
meeting the $x$-axis (resp. $y$-axis) with multiplicity at least $n-j+1$
(resp. $j+1$).\par
To see the same thing schematically, via equations in the local model
$H_{n+1}$, we proceed as follows. We work near
a \emph{generic} point $z_0\in D_{\theta'}$, necessarily of  multiplicity 1 at
the origin. Then we can, discarding distal
factors supported away from the nodes, write
the singleton (length-1) scheme corresponding to
a nearby cycle $z$ as  $\sch(z)=(x-c,y)$ where $c\to 0$ as $z\to z_0$, and
then
$$I_s(z)=(sx^{n-j}+y^j)(x-c,y)=(sx^{n-j+1}-csx^{n-j}-cy^j, y^{j+1}).$$ Thus, in
terms of the system of generators \eqref{Y-equations}
et seq., $I_s(z)$ is defined locally by \beql{I_s-equation}{cu_j-sv_j=0} (with
other $[u_k,v_k]$ coordinates either
$[1,0]$ for $k<j$ or $[0,1]$ for $k>j$.
The limit of this as $c\to 0$ is $[u_j,v_j]=[1,0]$, which is the point $Q_j$.
\emph{QED Claim.}\par
 Clearly $I_s$ doesn't meet $Q_j$ or $Q_{j+1}$ away from
 $D_{\theta'}\cup D_{\theta"}$. Therefore, denoting the
 scroll projection by $\pi$, we have
\beql{}{ I_s\cap Q_j=Q_j.\pi^*(D_{\theta'}),}
\beql{}{I_s\cap Q_{j+1}=Q_{j+1}.\pi^*(D_{\theta"});}
an easy calculation in the local model shows that the
intersection is transverse.
Because $Q_j\cap Q_{j+1}=\emptyset$, it follows that
\begin{eqnarray}I_a\sim Q_j+\pi^*(D_{\theta'})\\
I_a\sim Q_{j+1}+\pi^*(D_{\theta"}).\end{eqnarray}
These relations also follow from the fact,
which comes simply from setting $s=0$ or
dividing by $s$ and setting $s=\infty$ in \refp{I_s-equation}., that
\beql{limits}{ \lim\limits_{s\to 0}I_s=Q_j+\pi^*(D_{\theta'}),
\lim\limits_{s\to\infty} I_s= Q_{j+1}+\pi^*(D_{\theta"})}
It then  follows that
$$(Q_j)^2=Q_j.(I_s-\pi^*(D_{\theta'}))=Q_j.(Q_{j+1}+\pi^*(D_{\theta"}-D_{\theta'})),$$
hence
\begin{eqnarray}(Q_j)^2=Q_j.\pi^*(D_{\theta"}-D_{\theta'}),\end{eqnarray}
therefore finally
\beql{nodescroll-0}{L_j^n-L^n_{j+1}=\pi^*(D_{\theta"}-D_{\theta'}).}
This proves the Lemma.\end{proof}

 \newsection{Definition of node scrolls and
 polyscrolls}\label{globalization}
We now begin to extend our scope to a global proper
family $X/B$ of nodal curves,
 with possibly higher-dimensional base and fibres with
more than one node. Our main interest is in the node
scrolls in this generality,
where, rather than living over a symmetric product, they
 become  $\P^1$-bundles over a relative Hilbert
scheme (of lower degree)
associated to a 'boundary family' of $X/B$, i.e a
family obtained, essentially, as the partial normalization
of the subfamily of $X/B$ lying over the normalization of a component of
the locus of singular curves in $B$ (viz. the boundary of $B$).
For our purposes, it will be convenient to work 'node by node',
associating to each a boundary family.
We begin by making the appropriate notion of boundary family precise.
\subsection{Boundary data}
Let $\pi:X\to B$ now denote an arbitrary flat
 family of
nodal curves of arithmetic genus $g$ over an irreducible base, with
smooth
generic fibre. In order to specify the additional information
required to define a node scroll, we make
the following definition.
\begin{defn}
A \underline{boundary datum} for $X/B$ consists of\begin{enumerate}
\item an irreducible variety $T$ with a map  $\delta:T\to B$
unramified
to its image;
\item a 'relative node' over $T$, i.e.
 a lifting $\theta:T\to X$ of $\delta$ such that each
$\theta(t)$ is a node of $X_{\delta(t)}$;
\item a labelling, continuous in $t$, of the two branches
of $X_{\delta(t)}$ along $\theta(t)$ as $x$-axis and $y$-axis,
denoted $X', X"$.

\end{enumerate}
Given such a datum, the \underline{associated boundary family
} $X^\theta_T$ is the
normalization (= blowup) of the base-changed family $X\times_BT$
along the section $\theta$, i.e.
\[ X^\theta_T=\Bl_\theta(X\times_BT),
\] viewed as a family of curves of arithmetic genus
$g-1$ with two smooth, everywhere distinct, individually defined
 marked points $\theta_x, \theta_y$ on the
 respective branches $X', X"$. We denote by $\phi$ the
 natural map fitting in the diagram  $$
 \begin{matrix}
 X^\theta_T&&\\
 \downarrow&\stackrel{\phi}{\searrow}&\\
 X\times_BT&\to&X\\
 \downarrow&&\downarrow\\
 T&\stackrel{\delta}{\to}&B.
 \end{matrix}
 $$

\end{defn}
\begin{rem}
Note that the fibres of $X^\theta_T$ are disconneted
(e.g. a disjoint union of smooth curves of genera $i, g-i$) whenever $\theta$ is a separating node;
still they always have arithmetic genus $g-1$, where the arithmetic genus of a curve $X$
is defined as $1-\chi(\O_X)$.
\end{rem}
Note that a boundary datum indeed lives over the boundary of $B$;
in the other direction, we can associate to any
component $T_0$ of the boundary of $B$ a finite number
 of boundary data in this
sense: first consider a component $T_1$ of the
normalization of $T_0\times_B\mathrm{sing}(X/B)$, which already
admits a node-valued lifting $\theta_1$ to $X$, then further base-change
by the normal cone of $\theta_1(T_1)$ in $X$
(which is 2:1 unramified, possibly disconnected, over $T_1$),
to obtain a boundary datum as above. 'Typically', the curve corresponding
to a general
point in $T_0$ will have a single node $\theta$ and then the degree of
$\delta$ will be 1 or 2 depending on whether the branches along $\theta$
are distinguishable in $X$ or not (they always are distinguishable if
$\theta$ is a separating node
and the separated subcurves have different genera). Proceeding
in this way and taking all components which arise, we obtain finitely many
boundary data which 'cover', in an obvious sense,
the entire boundary of $B$. Such a collection, weighted so that
each boundary component $T_0$ has total weight $=1$ is called a
\emph{covering system of boundary data}.
\subsection{Node scrolls: definition}
\begin{propdef}
Given a boundary datum $(T, \delta, \theta)$ for $X/B$ and
natural numbers $1\leq j<n$, there exists a $\P^1$-bundle
$F^{m,n}_j(\theta)$,
called a \underline{node scroll} over
the Hilbert scheme $(X^\theta_T)\sbr m-n.$,
endowed with two disjoint sections $Q^{m,n}_{j,j}(\theta),
Q^{m,n}_{j+1,j}(\theta)$,
 together with a surjective
map generically of degree equal to $\deg(\delta)$
of
$$\bigcup\limits_{j=1}^{n-1} F^{m,n}_j(\theta):=
\coprod _{j=1}^{n-1} F^{m,n}_j(\theta)/
\coprod _{j=1}^{n-2} (Q^{m,n}_{j+1,j}(\theta)\sim Q^{m,n}_{j+1,
j+1}(\theta))$$
onto
the closure in $X\sbr m._B$ of the locus of
schemes having length precisely $n$ at $\theta$, so that
a general fibre of $F^{m,n}_j(\theta)$ corresponds to the family $C^n_j$
of length-$n$ schemes at
$\theta$ generically of type $C^n_j$, with the two nonprincipal schemes
$Q^{n}_j, Q^{n}_{j+1}$ corresponding to
$Q^{m,n}_{j,j}(\theta), Q^{m,n}_{j,j+1}(\theta)$ respectively.
We denote by $\delta^n_j$ the natural map of $F^{m,n}_j(\theta)$
to $X\sbr m._B$.
\end{propdef}\begin{proof}[Proof-construction]
We fix $m$ and $\theta$ (and suppress them when convenient). The scroll
$F_j^{m,n}(\theta)$
is defined as follows. Fixing
the boundary data, consider first the locus \[\bar{ F}_j^n \subset T\times_B
X\sbr m._B\]  consisting of
compatible pairs $(t,z)$ such that $z$ is in the closure of
the set of schemes which are of type $I^n_j$ (i.e. $x^{n-j}+ay^j, a\in \C^*$)
at $\theta(t)$, with
respect to the branch order $(\theta_x, \theta_y)$.
The discussion of \S\ref{study of H_m} shows that
the general fibre of $\bar F_j$ under the cycle map is a
$\P^1$, namely a copy of $C^n_j$; moreover the closure
of the locus of schemes having multiplicity $n$ at $\theta$
is the union $\bigcup\limits_{j-1}^{n-1}\bar F^{m,n}_j.$
In fact locally over a neighborhood of a
cycle having multiplicity precisely $n+e$ at $\theta$, $\bar F^{m,n}_j$
is a union of components $\bar F^{n}_j(a,b)\times U\spr m-e., a+b=e$,
where $U$ is an open set disjoint from $\theta$,
$\bar F^{n}_j(a,b)\subset H_{n+e}$ maps to $(U')^a\times (U")^b$ and is
defined in  $H_{n+e}$ by
 by the vanishing of
all $Z_i, i\neq j+b, j+b+1$ or alternatively, in terms of $u,v$ coordinates,
by
\[ v_1=...=v_{j+b}=u_{j+b+1}=...=u_{n+e}=0\]\par
Then $F^{m,n}_j(\theta)$ is the locus \eqspl{}{
\{(w,t,z)\in (X^\theta_T)\sbr m-n.\times_T\bar{F}^n_j:
\phi_*(c_{m-n}(w))+n\theta=c_m(z)\},
}
where $\phi:X^\theta\to X$ is the natural map, clutching together
$\theta_x$ and $\theta_y$, and $\phi_*$ is the induced push-forward map
on cycles. Then the results of the previous section show that
$F^{m,n}_j(\theta)$  is
locally defined near a cycle
having multiplicity $b$ at $\theta_y$, e.g. by the vanishing
of the $Z_i, i\neq j+b, j+b+1$ on
$$\{(w,u,Z)\in (X^\theta_T)\sbr m-n.\times X\spr
e._B\times\P^{n+e}:\phi_*(c_{m-n}(w))_\theta+n\theta=u\}$$
where $\underline{\ \ }_\theta$ indicates the portion near $\theta$.
The latter locus certainly projects isomorphically to its image
in $(X^\theta_T)\sbr m-n.\times\P^{n+e}$, hence $F^{m,n}_j(\theta)$ is a
$\P^1$-bundle over $(X^\theta_T)\sbr m-n.$.
Since $F^{m,n}_j(\theta)$
admits the two sections $Q^{m,n}_{j,j}(\theta), Q^{m,n}_{j+1,j}(\theta)$,
it is the projectivization of a decomposable rank-2
vector bundle.
\end{proof}
Note that the node scroll $F^{m,n}_j(\theta)$ also depends on $m$, and
is
by construction a subscheme of the 'flag-like' Hilbert scheme
\eqspl{nodescroll-corresp}{
\begin{matrix}
F^{m,n}_j(\theta)\subset&\{(z_1, z_2):\phi(z_1)\subset z_2\}
&\to  X\sbr m._B\\&\downarrow&\\
&(X^\theta_T)\sbr m-n.&
\end{matrix}}
Of course $z_1, z_2$ live on different families so this is not
the usual flag-Hilb.
We will denote the two Hilbert-scheme targeted projections
on $F^{m,n}_j(\theta)$ by $p_{[m-n]}, p_{[m]}$ respectively.
When the dependence on $\theta, m,...$ is obvious, we will
omit the corresponding designator.
The following simple technical point will be needed below.
\begin{lem}\label{disjoint-sections}Let
$T'\to T$ be a base change and $\theta'$ a
section
of $X^\theta_{T'}$ disjoint
from the distinguished sections $(\theta_x)_{T'}, (\theta_y)_{T'}$ and
identified with the corresponding section of
$X_{T'}$. Then on the
pulled- back node scroll $F^{m,n}_j(\theta)_{T'}$,
\[ p_{[m]}^*[m]_*\theta'=p_{[m-n]}^*[m-n]_*\theta'\]
\end{lem}
\begin{proof}
It suffices to verify this on the ordered version where, e.g.
$[m]_*\theta'=\sum\limits_{i=1}^m p_i^*\theta'$ and the projection
$p_{[m-n]}$ corresponds to projection on the first $m-n$
coordinates. But then for
$i>m-n$,  we have $p_i^*\theta'\cap F=\emptyset$
as the nodes are disjoint. This gives our assertion.
\end{proof}
Obviously, $Q^{m,n}_{j,j-1}(\theta)$ and $Q^{m,n}_{j,j}(\theta)$ coincide
in
$(X^\theta_T)\sbr m-n.\times X\sbr m._B$ and when convenient we
will write them as $Q^{m,n}_j(\theta)$ or $Q^{m,n}_j(\theta)$,
omitting $\theta$ when harmless. It is noteworthy
that the map from $Q^{m,n}_j(\theta)$ can be written down explicitly:
\begin{lem}\label{q-map}
The map  $(X_T^\theta)\sbr k.\simeq Q^{m,n}_j(\theta)\to X\sbr m._B$ is
given by
\eqspl{map-on-Q}{
 z_0+ a_x\theta_x+a_y\theta_y\mapsto \phi(z_0)+Q^{n+a_x+a_y}_{j+a_y}
 }
where $z_0$ is supported off $\theta_x\cup\theta_y$.
\end{lem}
\begin{proof}
To begin with, as $\theta_x, \theta_y$ are smooth sections
of $X^\theta_T$, any length-$k$ subscheme of it can indeed be
expressed
uniquely as in the formula.
The formula is clearly true
when $a_x=a_y=0$. Then the general case follows by taking limits,
in view of the explicit local description of the schemes of type
$Q^p_r$ as
$(x^{p-r+1}, y^r)$.
\end{proof}

\subsection{Polyscrolls}\label{polyscroll-sec}
Consider now a collection $\theta_.=(\theta_1,...,\theta_r)$
of distinct relative nodes of $X/B$ and $T=T(\theta_1,...,\theta_r)
\to B$ a common boundary locus
for them, compatible with the boundary data
for each $\theta_i$. Thus, $X_T$ is endowed with
 $r$ distinct relative nodes that we still denote by
 $\theta_1, ..., \theta_r$. Let $X^{\theta.}_T$ be the blowup or
partial normalization of $X_T$
in $\theta_1,...\theta_r$. As the $\theta_i$
are disjoint, the blowing up may be done inductively,
in any order, or simultaneously. Let $(j.),(n.)$ be
sequences of $r$ positive integers with $(j.)<(n.)$ in the
sense that $j_i<n_i,\forall i$.  We aim to define a node polyscroll
$F:=F^{m,n.}_{j.}(\theta.; X/B)$. This can be done using induction on $r$.
Assume the $(r-1)$- polyscroll
$F'=F^{m-n_1,n_2,...,n_r}_{j_2,...,j_r}(\theta_2,...,\theta_r;
X^{\theta_1}_{T(\theta_1)})$
is defined, together with  maps
\eqsp{
\begin{matrix}
&F'&\stackrel{p_{[m-n_1]}}{\to}&(X^{\theta_1})\sbr m-n_1._{T(\theta_1)}\\
p_{[m-|n.|]}&\downarrow&&\\
&(X_T^{\theta.})\sbr m-|{n.}|.&&
\end{matrix}
}
the horizontal one being generically finite
and the vertical one a $(\P^1)^{r-1}$-bundle projection.
Of course, the node scroll $F^{m,n_1}_{j_1}(\theta_1;X/B)$
is a $\P^1$-bundle over $(X^{\theta_1})\sbr m-n_1._{T(\theta_1)}$. Define
$F$
as the fibre product
\eqspl{polyscroll-def}{
\begin{matrix}
&&F&&\\
&\swarrow&&\searrow&\\
F^{m,n_1}_{j_1}(\theta_1;X/B)&&\diamondsuit&&F'\\
&\searrow&&\swarrow&\downarrow\\
&&{(X^{\theta_1})\sbr m-n_1._{T(\theta_1)}}&&(X_T^{\theta.})\sbr m-|{n.}|.
\end{matrix}
}
Then $F$ comes equipped with a $(\P^1)^r$-bundle projection
$p_{[m-|n.|]}F\to F'\to (X^{\theta.}_T)\sbr {m-{|n.|}}.$, as well as a
generically
finite map $p_{[m-n_1]}:F\to F^{m,n_1}_{j_1}(\theta_1;X/B) \to X\sbr
m._B$.
Writing, suggestively, $F'$ as $F(\hat\theta_1)$ and assuming
inductively
maps $F'\to F'(\hat\theta_i), \forall i>1$, we can identify
$F^{m,n_1}_{j_1}(\theta_1;X/B)\times F'(\hat\theta_i)$ as
$F(\hat\theta_i)$ and obtain an induced map
$F\to F(\hat\theta_i)$. Then taking fibre product with
$F^{m,n_i}_{j_i}(\theta_i;X/B)$, we obtain a morphism,
easily seen to be an isomorphism, from $F$ to a similar node polyscroll
with $\theta_1, \theta_i$ interchanged. Continuing
in this way,
 it is easy to see that $F$ is independent of
the ordering and the composite
$F\to F'\to (X^{\theta.}_T)\sbr {m-{|n.|}}.$
is a $(\P^1)^r$-bundle.
\par

We summarize some of the important properties of node polyscrolls as
follows
\begin{prop}\begin{enumerate}\item
The $r$-polyscroll $F=F^{m,n.}_{j.}(\theta..X/B)$ is a $(\P^1)^r$-bundle
over
the Hilbert scheme $(X_T^{\theta.})\sbr m-{|n.|}.$.
\item $F$ parametrizes subschemes of $X/B$ having length at least $n_i$
    at $\theta_i, i=1,...,r$.\item $F$ is independent of the order of $(\theta.,
n.,j.)$ and admits a $(\P^1)^{r-s}$-bundle projection to a pullback of the
$s$-polyscroll
based on any $s$ of the $(\theta_i, n_i, j_i)$.
\end{enumerate}
\end{prop}


 \newsection{Structure of node scrolls}\label{globalization2}
We fix a boundary datum $(T,\delta, \theta)$ as above.
Our aim now is to determine the structure of a node scroll as
$\P^1$ bundle together with the relative polarization induced by
minus the discriminant.
The following result is critical:
\begin{prop}\label{gamma-on-Q}
Let $Q^{m,n}_j=Q^{m,n}_j(\theta)$ be the canonical section of type
$Q^{m,n}_j$ in the node scroll
\[F^{m,n}_j=F^{m,n}_j(\theta)\subset X\sbr m._B.\] Then up to linear
equivalence, we have, where $k=m-n$, $Q^{m,n}_j$ is identified with
$(X^\theta)\sbr k._T$ and $\Gamma\spr k.$ is the discriminant on the
latter:
\eqspl{n general}{
\Gamma\spr m..Q^{m,n}_j\equiv_{\mathrm{lin}}
-\binom{n-j+1}{2}\psi_x-\binom{j}{2}
\psi_y+(n-j+1)[k]_*\theta_x+j[k]_*\theta_y+\Gamma\spr k.
:=D^{m,n}_j(\theta)}

\end{prop}
\begin{proof}
We begin with the \emph{special case $n=2$}. Here the possible values of
$j$ are 1 and 2 and by symmetry it suffices to consider $j=1$, where the
formula
reads
\eqspl{n=2}{\Gamma\spr m..Q^{m,2}_1\sim
-\psi_x+2[k]_*\theta_x+[k]_*\theta_y+\Gamma\spr k..}
Recall that $Q=Q^{m,2}_1$ is the graph of the morphism
$q: (X^\theta)\sbr k._T\to X\sbr m._B$ given as in Lemma \ref{q-map}.
Every scheme in the image of $q$ contains the length-2 scheme
along the $x$-axis, $(2\theta_x)$, locally defined by $(y,x^2)$.
Therefore $q$ clearly factors through a map
\[q':(X^\theta)\sbr k._T\to X\sbr 2,m._B\]
to the Hilbert scheme of $(2,m)$-flags.
Moreover the projection of $q'$ to $X\sbr 2._B$ is the relatively constant
map with value $(2\theta_x)$ (the unique length-2
scheme contained in the $x$-branch).\par
Now, $X\sbr 2,m._B$ carries the pullback of $-\Gamma\spr 2.$ from
$X\spr 2._B$, which
clearly pulls back via $q'$ to the cotangent space in the $x$ direction,
i.e.  $\psi_x$. So we get an injection
\[\O_{Q^{m,2}_1}(-\Gamma\spr m.)\subset
\psi_x\otimes\O_{Q^{m,2}_1}(-\Gamma\spr k.)\]
(where $\Gamma\spr k.=\Gamma_{X^\theta_T}\spr k.$ denotes
as above the discriminant associated
to the boundary family $X^\theta_T$).
This injection is clearly an iso over the open set of subschemes
of $X^\theta$ disjoint form $\theta_x\cup\theta_y$. Therefore
\[\O_{Q^{m,2}_1}(-\Gamma\spr m.) =
\psi_x\otimes\O_{Q^{m,2}_1}(-\Gamma\spr k.-\alpha
[k]_*\theta_x-\beta[k]_*\theta_y )\]
for some nonnegative integers $\alpha,\beta$. To identify these, we can
work at a general point of their support, which corresponds to a
scheme with a length-3 portion near $\theta$. By the usual support
decomposition of Hilbert schemes as in \S\ref{preliminary_reductions}, we
are reduced to the
case $m=3$, working near a scheme of type $Q^3_2=(y^2, x^3)$ for
$\beta$
 or
$Q^3_1=(y^1, x^2)$ for $\alpha$. Moreover, pulling back by the
 finite flat morphism from the ordered Hilb $X\scl 3._B$, we are reduced
 to working there with the 3rd coordinate being the one
 from $X^\theta$ and the first two corresponding to the map to
 $X\sbr 2._B$ (so that $y_1=y_2=x_1^2=x_2^2=0)$. \par
 Then finally, in the first case, the generator $G^3_2$
(i.e. the mixed Van der Monde)
 can be expanded along the last row, which shows that
 it maps to $y\psi_x$, therefore $\beta=1$. In the second case, the
 generator $G^3_1$ maps to
$x^2\psi_x$, so $\alpha=2$. This completes the proof in the case
$n=2$.\par Passing to the general case, recall from \S\ref{study of H_m}
that $Q^{m,n}_j$ is the pullback of $Q^{2,m}_1$ via the map
\eqsp{f:(X^\theta)\sbr k._B\to (X^\theta)\sbr m-2._B\\
z\mapsto z+(n-j-1)\theta_x+(j-1)\theta_y
}Then given \eqref{n=2}, the desired formula \eqref{n general} is
a consequence of following  elementary formulas (recall $k=m-n$)
\eqsp{
f^*([m-2]\theta_x)&=[k]_*\theta_x-(n-j-1)\psi_x\\
f^*([m-2]_*\theta_y)&=[k]_*\theta_y-(j-1)\psi_y\\
f^*(\Gamma\spr m-2.)&=\Gamma\spr k.+(n-j-1)[k]_*\theta_x
-\binom{n-j-1}{2}\psi_x+(j-1)[k]_*\theta_y-\binom{j-1}{2}\psi_y.
} Note that because $\theta_x, \theta_y$ map to a node of $X/B$, 
they are contained in the smooth part of $X^\theta/B$.
Then, note that $f$ is an iterate of maps of the following form, associated
to a section $\sigma:B\to Y$ of a nodal family 
\eqsp{i_\sigma:Y\sbr k._B\to Y\sbr k+1._B\\ i_\sigma(z)=z+\sigma
} To prove the above formulas, it suffices by an evident
recursion to prove the following Lemma,
which will conclude the proof of Lemma \ref{gamma-on-Q}.\end{proof}
\begin{lem}\label{sublemma}
For $\sigma$ as above and a section $\sigma'$ disjoint from $\sigma$,
assume $\sigma, \sigma'$ are contained in the smooth part of $Y/B$. Then
we have, where $\psi_\sigma=\omega_{Y/B}|_\sigma$
\eqspl{}{
i_\sigma^*([k+1]_*\sigma)=[k]_*\sigma-\psi_\sigma,\quad
i_\sigma^*([k+1]_*\sigma')=
[k]_*\sigma'\\
i_\sigma^*(\Gamma\spr k+1.)=\Gamma\spr k.+[k]_*\sigma
}
\end{lem}
\begin{proof}[Proof of \ref{sublemma}]
It suffices to prove the analogous fact on the relative symmetric product,
where = becomes linear equivalence of Weil divisors.
Because such linear equivalence descends via a finite flat map like the
symmetrization, it suffices to prove the analogous fact over the
relative Cartesian product. There, the first two assertions are
obvious (keeping in mind the the image of our sections is disjoint
from the nodes). The last assertion becomes obvious as well once
we recall that the big diagonal on the Cartesian product is the sum
of pullbacks from 2-fold products.
\end{proof}

We are now in position to determine the polarized structure of the
node scroll $F^{m,n}_j(\theta)$. This means finding a vector bundle $E$
such that $F^{m,n}_j=\P(E)$
and such that the canonical $\O(1)$ polarization on $\P(E)$
corresponds to $-\Gamma\spr m.$. We recall (see EGA or \cite{hart}, Ch.
II.7 or \cite{ful} which unfortunately uses the opposite
sign convention) that for any vector bundle $E$, there is a canonically
defined
(depending on $E$) line bundle on $\P(E)$, denoted $\O(1)$ or
$\O_E(1)$, which restricts to the usual (Grothendieck, or quotient) $\O(1)$ on (geometric)
fibres.
\begin{thm}[Node scroll theorem]\label{nodescroll-thm}
For any boundary datum $(T,\delta, \theta)$,
and any \mbox{$1\leq j<n\leq m$,} there is an isomorphism
\eqspl{polscrollsctruct}{
F^{m,n}_j(\theta)\simeq \P(\O(D^n_j(\theta))\oplus \O(D^n_{j+1}(\theta)))
}which pulls back the canonical $\O(1)$ polarization on the RHS
 to the restriction of $-p_{X\sbr m._B}^*\Gamma\spr m.+
p_{(X_T^{\theta.})\sbr m-n.}^*\Gamma\spr m-n.$ on the LHS.
\end{thm}
\begin{proof}
As $F^{m,n}_j(\theta)$ admits the two disjoint sections
$Q^{m,n}_j, Q^{m,n}_{j+1}$, the result is immediate from
Proposition \ref{gamma-on-Q}.
\end{proof}
\begin{cor}
On $F^{m,n}_j(\theta)$, we have
\eqspl{}{
-\Gamma\spr m.&\sim Q^{m,n}_j+p_{[m-n]}^*(D^n_{j+1})\\
&\sim Q^{m,n}_{j+1}+p_{[m-n]}^*(D^n_j).
}
\end{cor}
\begin{proof}
Follows from the elementary fact that on any $\P^1$-bundle $\P(A\oplus
B)$ with projection $\pi$, we have
\[c_1(\O(1))\sim \P(A)+\pi^*(c_1(B)).\]
Indeed the natural map $\pi^*(B)\to\O(1)$ vanishes precisely on the divisor
$\P(A)\subset\nobreak\P(A\oplus\nobreak B)$.
\end{proof}
The extension to polyscrolls is direct from the
definition in \S\ref{polyscroll-sec} once we note that thanks to the
disjointness
of the nodes, the divisors
$D^{m,n_i}_{j_i}(\theta_i)$ correspond naturally to a similarly-denoted
divisor on $(X_T^{\theta.})\sbr m-{|n.|}.$, with $[m-n_i]_*\theta_{x,i}$
corresponding to $[m-|n.|]_*\theta_{x,i}$, e.g. on
$F^{n_1}_{j_1}(\theta_1)$,
\[p_{[m]}^*[m]_*\theta_{2,x}=p_{[m-n_1]}^*[m-n_1]_*\theta_{2,x}\]
etc. where $p_{[k]}$ denotes the natural map to the
length-$k$ Hilbert scheme (of $X$ or $X^{\theta_1}$)
(compare Lemma \ref{disjoint-sections}).

\begin{thm}[Node Polyscroll Theorem]\label{nodepolyscroll-thm}
There is an isomorphism
\eqspl{}{
F^{m,n.}_{j.}(\theta.;X/B)\sim {\prod}_{(X^{\theta.}_{T(\theta.)})\sbr{m-|n.|}.\ \ } \P(\O(D^{m,n_i}_{j_i}(\theta_i))
\oplus \O(D^{m,n_i}_{j_i+1}(\theta_i)))
} under which $-\Gamma\spr m.+\Gamma\spr m-{|n.|}.$ corresponds
to the canonical $\O(1,...,1)$.
\end{thm}
\begin{proof}
We use the setting and notations of \S\ref{disjoint-sections}.
Consider the natural  projection
$F\to F'$, which is just a base-change of the scroll
projection \[p_{[m-n_1]}: F_1=F^{m,n_1}_{j_1}(\theta_1, X/B)
\to (X^{\theta_1})\sbr m-n_1..
\]
 Via this, we have
\eqsp{\O_{F_1}(1)=\binom{n_1-j_1+1}{2}\psi_{1,x}+\binom{j_1}{2}\psi_{1,y}
-(n_1-j_1+1)[m-n_1]_*(\theta_{1,x})-j_1[m-n_1]_*\theta_{1,y}.}
On $F$ this becomes, using Lemma \ref{disjoint-sections} (essentially,
the disjointness of the sections $\theta_i$),
\eqsp{
&\O_{F_1}(1)|_F=\\
&\binom{n_1-j_1+1}{2}\psi_{1,x}+\binom{j_1}{2}\psi_{1,y}
-(n_1-j_1+1)[m-|n.|]_*(\theta_{1,x})-j_1[m-|n.|]_*(\theta_{1,y}).}
and by Theorem \ref{nodescroll-thm}, this coincides on $F$ with
$-\Gamma\spr m.+\Gamma\spr m-n_1.|_{F'}$. By induction,\nl
$-\Gamma\spr m-n_1.|_{F'}+\nobreak\Gamma\spr{m-|n.|}.$ coincides
with
the appropriate $\O(1,...,1)$ on the $(r-1)$-polyscroll $F'$, and the
Theorem follows.

\end{proof}
\begin{rem}
As mentioned in the Introduction, the paper \cite{internodal} and the related software Macnodal
\cite{macnodal}
contain numerous numerical examples and applications of the Node Scroll and Polyscroll Theorems.
\end{rem}
\begin{rem}
Define a \emph{smudgy curve} of type $g,p,k$ to be a nodal, $p$-pointed,
genus-$g$ curve together with a length-$k$ subscheme such that the
entire object has finite automorphism
group, and
let $\overline\M\sbr k._{g,p}$ denote the moduli space
(DM stack) of smudgy curves of this type (assuming it exists). Some
interesting questions about (ordinary) curves (for example, Brill-Noether
loci) can be formulated in terms of smudgy curves. The node scrolls define
correspondences between smudgy moduli spaces:
\[\overline\M\sbr k_1._{g_1,p_1+1}\times \overline\M\sbr
k_2._{g_2,p_2+1}
{\leftarrow} \ \ \ F^{m,n}_j\ \ \to\overline\M\sbr
k_1+k_2+n._{g_1+g_2,p_1+p_2},\]
\[\overline\M\sbr k._{g-1,p+2}
{\leftarrow} \ \ \ F^{m,n}_j\ \ \to\overline\M\sbr k+n._{g,p}\]
where $k=k_1+k_2, p=p_1+p_2, g=g_1+g_2$
(identifying the LHS with a boundary component of $\overline\M\sbr
k._{g,p}$).  These are analogous to the correspondences used by
Nakajima \cite{nak} to define creation-annihilation operators on the
cohomology of Hilbert schemes of surfaces.
\end{rem}
\vfill\eject
\pagestyle{plain}
\bibliographystyle{amsplain}
\bibliography{../mybib}

\providecommand{\bysame}{\leavevmode\hbox to3em{\hrulefill}\thinspace}
\providecommand{\MR}{\relax\ifhmode\unskip\space\fi MR }
\providecommand{\MRhref}[2]{%
  \href{http://www.ams.org/mathscinet-getitem?mr=#1}{#2}
}
\providecommand{\href}[2]{#2}
\begin{thebibliography}{10}

\bibitem{ang}
B.~Ang\'eniol, \emph{Familles de cycles alg\'ebriques- sch\'ema de {Chow}},
  LNM, vol. 896, Springer.

\bibitem{acgh}
E.~Arbarello, M.~Cornalba, P.~Griffiths, and J.~Harris, \emph{Geometry of
  algebraic curves}, Grundl. d. math. Wiss., vol. 267, Springer, 1985.

\bibitem{cotteril1}
E.~Cotteril, \emph{Geometry of curves with exceptional secant planes}, Math. Z.
  \textbf{267} (2011), 549--582.

\bibitem{ful}
W.~Fulton, \emph{Intersection theory}, Ergeb. d. Math. u. i. Grenzgeb. 3.
  Folge, Bd. 2, Springer, Berlin, 1984.

\bibitem{hart}
R.~Hartshorne, \emph{Algebraic {G}eometry}, Springer, Berlin, 1977.

\bibitem{kl}
S.~Kleiman and D.~Laksov, \emph{On the existence of special divisors}, Amer. J.
  Math (1972), 431--436.

\bibitem{L}
M.~Lehn, \emph{{Chern} classes of tautological sheaves on {Hilbert} schemes of
  points on surfaces}, Invent. math \textbf{136} (1999), 157--207.

\bibitem{lehn-montreal}
\bysame, \emph{Lectures on {H}ilbert schemes}, CRM notes, Centre de
  {R}echerches {M}ath\'ematiques, Montreal, 2004.

\bibitem{macnodal}
Gwoho Liu, \emph{The macnodal package for intersection theory on hilbert
  schemes of nodal curves}, web interface (small jobs only) at
  \url{http://gwoho.com/macnodal/index.html.a?format=html}; source + executable
  + instructions at \url{http://math.ucr.edu/\~ziv/}.

\bibitem{macd}
I.~G. MacDonald, \emph{Symmetric products of an algebraic curve}, Topology
  \textbf{1} (1962), 319--343.

\bibitem{nak}
H.~Nakajima, \emph{Lectures on {H}ilbert schemes of points on surfaces},
  University lecture series, Amer. Math. Soc., 1999.

\bibitem{grd}
Z.~Ran, \emph{Boundary modifications of {H}odge bundles and enumerative
  geometry}, \url{http:arXiv.org/abs/1011.0406v2}.

\bibitem{R2}
\bysame, \emph{Cycle map on {Hilbert} schemes of nodal curves}, Projective
  varieties with unexpected properties (Ciliberto et~al., ed.), De Gruyter,
  Berlin, 2005, pp.~363--380.

\bibitem{geonodal}
\bysame, \emph{Geometry on nodal curves}, Compositio math \textbf{141} (2005),
  1191--1212.

\bibitem{Hilb}
\bysame, \emph{A note on {Hilbert} schemes of nodal curves}, J. Algebra
  \textbf{292} (2005), 429--446.

\bibitem{internodal}
\bysame, \emph{Tautological module and intersection theory on {H}ilbert schemes
  of nodal curves}, Asian J. Math. \textbf{17} (2013), 193--264,
  \url{arxiv:0905.2229v5}.

\bibitem{sernesi}
E.~Sernesi, \emph{Deformations of algebraic schemes}, Grundl. d. math. Wiss.,
  vol. 334, Springer International, Berlin, Heidelberg, 2006.

\bibitem{Vak}
R.~Vakil, \emph{The moduli space of curves and {Gromov-Witten} theory},
  Enumerative invariants in algebraic geometry and string theory (K.~Behrend
  and M.~Manetti, eds.), Springer-Verlag, 2008,
  \url{arxiv.org/math.AG/0602347}.

\end{thebibliography}
\end{document}